\documentclass[12pt]{amsart}
\usepackage{latexsym}
\usepackage[all]{xy}
\usepackage{amsfonts}
\usepackage{amsthm}
\usepackage{amsmath}
\usepackage{amssymb}
\usepackage{amscd}
 \numberwithin{equation}{section}
 \usepackage{mathptmx}

\def\sw#1{{\sb{(#1)}}}
    \def\sco#1{{\sb{[ #1]}}} 
    \def\su#1{{\sp{[#1]}}} 
    
    \def\zz#1{{\sb{[#1]}}}
    \def\tens{\mathop{\otimes}}

    \def\<{{\langle}}
    \def\>{{\rangle}}
    
    \def\eps{\varepsilon}

    \def\note#1{{}}

    \def\note#1{}
    \def\M{{\bf M}}

    \def\cA{{\mathfrak A}}

    \def\eA{{\eta_\cA}}
    \def\mA{{\mu_\cA}}
    \def\roA{{\varrho^\cA}}
    \def\eB{{\eta_\cB}}
    \def\mB{{\mu_\cB}}
    \def\Bro{{}^\cB\!\varrho}
    \def\roB{{\varrho^\cB}}
    \def\Aro{{}^\cA\!\varrho}

    \def\cA{{\mathcal A}}
     \def\cX{{\mathcal X}}

    \def\beq{\begin{equation}}
    \def\eeq{\end{equation}}

    \def\im{{\rm Im}}

    \def\ot{{\otimes}}

    \def\roM{\varrho^{M}}
    \def\Nro{{}^{N}\!\varrho}
    \def\roN{\varrho^{N}}
    \def\Mro{{}^{M}\!\varrho}
     \def\roC{\varrho_{C}}

    \newcommand{\Ra}{\Rightarrow}

    \def\vec{\mathbf{Vect}_k}
    \def\wra{\mathop{\overline{\rho_M}}}
    \def\wrac{\mathop{\overline{\rho_C}}}
    \def\wla{\mathop{\overline{{}_N\rho}}}
    \def\wlac{\mathop{\overline{{}_C\rho}}}
    
    \def\lpr{\overline{p_R}}
    \def\lpl{\overline{p_L}}

    \def\sba{_{\alpha}}
    \def\spa{^{\alpha}}
    \def\sbb{_{\beta}}
    \def\spb{^{\beta}}
    \def\spab{^{\alpha\beta}}
    
    \def\sbab{_{\alpha\beta}}

    \def\sumab{\sum_{\alpha,\beta}}
    \def\suma{\sum_{\alpha}}
    \def\sumAB{\sum_{E,F}}
    \def\sumA{\sum_{E}}
    \def\cuc{\varepsilon_{C}}
    \def\cmc{\Delta_{C}}
    \def\re{\psi_{R}}
    \def\le{\psi_{L}}
    \def\rem{\M(\psi)_{A}^{C}}
    \def\lem{{_{A}^{C}\M(\le)}}
    \def\rema{\M(\re)_{A}^{C}}
    \def\htau{\widehat{\tau}}
    
    \def\cB{{\mathcal B}}
    \def\cten#1{\raise-.2cm\hbox{$\stackrel{\displaystyle\square}
{\scriptscriptstyle{#1}}$}}

\def\coten#1{\,\underset{\scriptscriptstyle{#1}}{\square}\,}
\def\cotc{\coten C}
    \def\cotb{\coten B}
    \def\cotd{\coten D}
    \def\cots{\coten{S}}
    \def\epi{\twoheadrightarrow}

     \newcounter{zlist}
  \newenvironment{zlist}{\begin{list}{(\arabic{zlist})}{
  \usecounter{zlist}\leftmargin2.5em\labelwidth2em\labelsep0.5em
  \topsep0.6ex
  \parsep0.3ex plus0.2ex minus0.1ex}}{\end{list}}

  \newcounter{blist}
  \newenvironment{blist}{\begin{list}{(\alph{blist})}{
  \usecounter{blist}\leftmargin2.5em\labelwidth2em\labelsep0.5em
  \topsep0.6ex 
  \parsep0.3ex plus0.2ex minus0.1ex}}{\end{list}}

\newtheorem{proposition}{Proposition}[section]
\newtheorem{lemma}[proposition]{Lemma}

\newtheorem{corollary}[proposition]{Corollary}
\newtheorem{theorem}[proposition]{Theorem}
\theoremstyle{definition}
\newtheorem{definition}[proposition]{Definition}
\newtheorem{example}[proposition]{Example}
\theoremstyle{remark}
\newtheorem{remark}[proposition]{Remark}

\begin{document}
\title{The Galois theory of matrix $C$-rings}
   \author{Tomasz Brzezi\'nski}
   \author{Ryan B.\ Turner}
   \address{ Department of Mathematics, University of Wales Swansea,
   Singleton Park, \newline\indent  Swansea SA2 8PP, U.K.}
   \email{ T.Brzezinski@swansea.ac.uk (T Brzezi\'nski)}

   \email{ marbt@swansea.ac.uk (RB Turner)}
     \date{\today}
   \begin{abstract}
 A theory of monoids in the category of bicomodules of a coalgebra $C$ or
 $C$-rings is developed. This can be viewed as a dual version of the coring theory.
 The notion of a matrix ring context consisting of two bicomodules and two maps
 is introduced and the corresponding example of a $C$-ring (termed a {\em matrix
 $C$-ring}) is constructed. It is shown that a matrix ring context can be associated to
 any bicomodule which is a one-sided quasi-finite injector. Based on this, the notion
 of a {\em Galois module} is introduced and the structure theorem, generalising
 Schneider's Theorem II [H.-J.\ Schneider,
Israel J.\ Math., 72 (1990), 167--195], is proven. This is then applied to the $C$-ring
 associated to a weak entwining structure and a structure theorem for a
 weak $A$-Galois coextension is derived. The theory of matrix ring contexts for a firm coalgebra
(or {\em infinite matrix ring contexts}) is outlined. A Galois
connection associated to a matrix $C$-ring is constructed.

   \end{abstract}
   \maketitle

\section{Introduction}
 The present paper is a contribution to the long standing programme (motivated
by non-commutative geometry) of understanding the origins and
finding the most general formulation of Schneider's structure
theorems for Galois-type extensions \cite{Sch:pri}. With the
re-birth of interest in corings triggered by \cite{Brz:str} it has become clear that the
proper general formulation of Schneider's Theorem I can be
provided by corings and their comodules, and such formulations
were achieved in recent papers \cite{Brz:gal}, \cite{CaeDeG:com},
\cite{Wis:gal}, \cite{Boh:Gal}. It had earlier been realised in
\cite{Brz:str} that to obtain a generalisation of Schneider's
Theorem II, which can be understood as a dual version of
Theorem~I, one needs to develop new algebraic structures, termed
{\em $C$-rings} in \cite[Section~6]{Brz:str}. Given a coalgebra
$C$ (over a field $k$), a {\em $C$-ring} is a monoid in the
category of $C$-bicomodules (with the monoidal structure provided
by the cotensor product $-\coten C -$). Explicitly a $C$-ring is a
$C$-bicomodule $\cA$ together with two bicomodule maps $\mA:
\cA\coten C\cA\to\cA$ and $\eA:C\to \cA$ such that
$$
\mA\circ(\mA\coten C\cA) = \mA\circ(\cA\coten C\mA), \qquad \mA\circ (\cA\coten C \eA) = \mA\circ(\eA\coten C\cA) = \cA,
$$
where  the standard isomorphisms $\cA\coten CC\simeq \cA \simeq
C\coten C\cA$ provided by the $C$-coactions are implicitly used.
The current most general formulation of Schneider's Theorem I
involves not so much corings themselves but a special class of
their comodules, termed {\em principal comodules}. Crucial for
this formulation is  the notion of a {\em comatrix coring}
introduced in \cite[Proposition~2.1]{ElKGom:com}, i.e.\ a coring
which can be associated to any bimodule, finitely generated and
projective on one side. Prompted by this in the present paper we
introduce and study {\em matrix $C$-rings}, which can be
associated to any $(D,C)$-bicomodule that is a quasi-finite
injector as a $C$-comodule. As the notion of a {\em quasi-finite
injector} is not as familiar as the notion of a finitely generated
projective module, in our definition of a matrix $C$-ring we
follow the route suggested by \cite[Theorem~2.4]{BrzGom:com}, and
define matrix $C$-rings through {\em matrix ring contexts}. The
latter have a very natural meaning as adjoint pairs in a
bicategory of bicomodules and are very closely related to
Morita-Takeuchi contexts \cite{Tak:Mor}.

Recall from \cite[Section~6]{Brz:str} that a right module of a $C$-ring $\cA$
is a right $C$-comodule $M$ together with a right $C$-comodule map
$\wra: M\coten C \cA\to M$ such that
$$
\wra\circ(\wra\coten C\cA) = \wra\circ(M\coten C\mA), \qquad
\wra\circ(M\coten C\eA) = M,
$$
where again the standard isomorphism
$M\coten CC\simeq M$ provided by the $C$-coaction on $M$
is implicitly used. We introduce the notion of
 an {\em $\cA$-coendomorphism coalgebra}
of a right $\cA$-module, quasi-finite and injective as a
$C$-comodule. Starting with a right $\cA$-module $M$ which is a
quasi-finite injector as a $C$-comodule, we are able to construct
a matrix $C$-ring. If this $C$-ring is isomorphic to $\cA$, then
we say that $M$ is a {\em Galois module}. If, furthermore, $M$ is
an injective module of the $\cA$-coendomorphism coalgebra, then we
say that $M$ is a {\em principal Galois module}. We then derive the
equivalent conditions for $M$ to be a Galois and principal Galois module
in Theorem~\ref{theorem1}. This is the main result of the paper,
and is a sought generalisation of Schneider's Theorem II. We then
construct a Galois connection associated to a matrix $C$-ring.
Finally we apply Theorem~\ref{theorem1} to a $C$-ring associated
to a weak entwining structure and obtain a dual version of results
in \cite{BrzTur:str}. In particular we prove that, within an
invertible weak entwining structure, a coextension of coalgebras
by a (left) self-injective algebra has a Galois property, provided
the canonical map is injective.\bigskip

\noindent{\bf Notation.} We work over a field $k$. For a coalgebra
$C$, the product is denoted by $\Delta_C$ and the counit by
$\eps_C$. For a right (resp.\ left) $C$-comodule $M$ the coaction
is denoted by $\roM$ (resp.\ $\Mro$). We use Sweedler's notation
for coproducts $\Delta_C(c) =c\sw 1\ot c\sw 2$, for right
coactions $\roM(m) = m\zz 0\ot m\zz 1$, and for left coactions
$\Mro (m) = m \zz{-1}\ot m\zz 0$. The cotensor product is denoted
by $-\cotc-$. For a $C$-ring $\cA$, $\mA$ denotes the product (as
a map, on elements it is denoted by a juxtaposition), $\eA$ is the
unit, $\wra$ (resp.\ $\overline{{}_M\varrho}$) is the $\cA$-action
on right (resp.\ left) $\cA$-module $M$. The categories of right
(resp.\ left) $\cA$-modules and $C$-comodules are denoted by
$\M_\cA$ and $\M^C$ (resp.\ ${}_\cA\M$ and ${}^C\M$), while
${}^C\M_A$ denotes the category of right $A$-modules and left
$C$-comodules with right $A$-linear coaction.

\section{Matrix ring contexts}
\subsection{Quasi-finite matrix contexts.}
\begin{definition} \label{MRC}
A {\em matrix ring context}, $(C, D, {^C N^D}, ^D\! {M^C}, \sigma, \tau)$,
consists of a pair of coalgebras $C$ and $D$, a $(C,D)$-bicomodule
 $N$, a $(D,C)$-bicomodule $M$ and a pair of bicomodule maps
$$
\sigma: C \to N \coten D M, \qquad \tau: M \coten C N \to D
$$
such that the diagrams

$$
\xymatrix{
N \coten D M \coten C N  \ar[d]_{N\coten D\tau} & &   C \coten C N
\ar[ll]_{\sigma\coten C N} &&   M \coten C N \coten D M   \ar[rr]^{\tau \coten D M}
&  &       D \coten D M\\
N\coten DD && N\ar[ll]_{\roN} \ar[u]_{\Nro} & & M\coten C C\ar[u]^{M\coten C
\sigma} & &M \ar[ll]_{\roM} \ar[u]_{\Mro}}
$$commute. The map $\sigma$ is called a {\em unit} and $\tau$ is called
a {\em counit} of a matrix context.
\end{definition}
Since a counit $\tau$ of a matrix ring context is a $D$-bicomodule map, it
is fully determined by its {\em reduced form} $\widehat{\tau} = \eps_D\circ\tau$.
The map $\widehat{\tau}$ is called a {\em reduced counit} of a matrix context.
Note that the $D$-bicolinearity of $\tau$ is equivalent to the following property
of $\htau$,
\begin{equation}\label{eq.htau.bicol}
(D\ot \htau)\circ(\Mro\ot N) = (\htau\ot D)\circ (M\ot \roN).
\end{equation}
In terms of the reduced counit, the commutative diagrams in Definition~\ref{MRC}
read
\begin{equation}\label{eq.MRC}
(N\ot \htau)\circ(\sigma\coten C N)\circ \Nro = N,
\qquad (\htau\ot M)\circ (M\coten C\sigma)\circ \roM = M.
\end{equation}
In other words, equations \eqref{eq.MRC} mean that $(N\ot \htau)\circ(\sigma\coten C N)$ is the identity on $C\coten C N$, while
$(\htau\ot M)\circ (M\coten C\sigma)$ is the identity on $M\coten C C$.

The notion of a matrix context is closely related to that of {\em pre-equivalence data}
or a {\em Morita-Takeuchi context} introduced in \cite[Definition~2.3]{Tak:Mor}.
In particular, in view of \cite[Theorem~2.5]{Tak:Mor}, if one of the maps
in a Morita-Takeuchi context is injective, then there is a corresponding
matrix ring context. Furthermore, every equivalence data  give rise to a matrix
ring context. This relationship explains the use of term {\em context} in
Definition~\ref{MRC}. The use of term {\em ring} is justified by the following
\begin{proposition}\label{prop.matring}
Let $(C, D, {^C N^D}, ^D\! {M^C}, \sigma, \tau)$ be a matrix ring context. Then
$\cA := N\coten DM$ is a $C$-ring with the product and unit
$$
\mA  = N \coten D\htau \coten D M, \qquad \eA = \sigma,
$$
where $\htau$ is the reduced counit. Furthermore, $M$ is a right $\cA$-module with the action $\htau\coten{D}M$ and $N$ is a left $\cA$-module with the action
$N\coten{D}\htau$. The $C$-ring $\cA$ is called a {\em matrix
$C$-ring}.
\end{proposition}
\begin{proof}
By definition, both $\mA$ and $\eA$ are $C$-bicomodule maps. Since
$\tau$ is a $D$-bicomodule map (cf.\ equation
\eqref{eq.htau.bicol}), the product $\mA$ is well-defined, i.e.\
$\mA (\cA\coten C \cA) \subseteq \cA$. The associativity of the
product $\mA$ follows immediately by the $k$-linearity of $\htau$,
while equations \eqref{eq.MRC} imply that $\eA = \sigma$ is the
unit for $\mA$. The statements about the actions of $\cA$ are
proven in a similar way.
\end{proof}

\begin{example}\label{ex.MRC.trivial}
As an immediate example of a matrix ring context, consider a coalgebra
map $f: C\to D$. Take $M=N=C$, viewed as a $(D,C)$- or $(C,D)$-bicomodule
via the map $f$, and define $\sigma =\Delta_C$ and $\tau = f$. Note that
$\htau = \eps_D\circ f = \eps_C$. The
corresponding matrix $C$-ring is $\cA = C\coten D C$ with the product
$\mu_A = C\coten D\eps_C\coten D C$ and
unit $\Delta_C$.
\end{example}

We now explore the meaning of a matrix context.

\begin{proposition} \label{mrcaf} If $(C, D, {^C N^D}, ^D {M^C}, \sigma, \tau)$
is a matrix ring context, then the cotensor
functor $F = {- \coten C} N : \M ^C \to \vec$ is a left adoint of the tensor functor $G = -
\otimes M: \vec \to \M ^C$.
\end{proposition}
\begin{proof}
Define natural transformations
$$
\varphi : \M^C \to GF, \qquad \varphi_X:= (X \coten C \sigma)
\circ \varrho^X,
$$
$$\nu : FG \to \vec, \qquad \nu_Y := (Y
\otimes \eps_D) \circ (Y \otimes \tau) = Y\otimes \widehat{\tau}.
$$
We need to show that these morphisms are the unit and counit,
respectively, of the adjunction. Take any right $C$-comodule $X$ and
compute
\begin{eqnarray*}
\nu_{F(X)} \circ F( \varphi_X)&=& (X\coten C N \otimes \htau)
\circ (((X \coten C \sigma) \circ \varrho^X )\coten C
N)\\
&=&(X\coten C N \otimes \htau) \circ
(X \coten C  \sigma \coten C N)\circ (X \coten C { \Nro})
= X \coten C N = F(X),
\end{eqnarray*}
where the  second equality follows by the definition of the cotensor product, and the third equality follows by the  first of equations \eqref{eq.MRC}.  On the other hand,
for all vector spaces $Y$,
$$
G(\nu_Y) \circ \varphi_{G(Y)} =
(Y \otimes \htau \otimes M) \circ (Y \otimes M \coten C \sigma)
\circ (Y \otimes \roM)
= Y \otimes M = {G(Y)},
$$
by the second of equations  \eqref{eq.MRC}.
Hence the natural transformations $\varphi$ and $\nu$ satisfy the
required properties.
\end{proof}

Since, given a matrix ring context $(C, D, {^C N^D}, ^D {M^C}, \sigma, \tau)$,
the functor $-\ot M: \vec\to \M^C$ has a left adjoint, the right $C$-comodule $M$ is
a {\em quasi-finite} comodule (cf.\ \cite[Proposition~1.3]{Tak:Mor}).
The left adjoint of $-\ot M:\vec\to \M^C$ is known as a {\em co-hom} functor and
is denoted by $h_C(M,-):  \M^C\to \vec$. By the uniqueness of adjoints, in
the case of a matrix ring context,
$h_C(M,-) \simeq -\coten C N$. Since $h_C(M,-)$ has a right adjoint,
it is right exact, and since $-\coten C N$ is left exact, the above isomorphism
of functors implies that the  cohom functor $h_C(M,-)$ is exact, i.e.\ the
right $C$-module
$M$ is an {\em injector} (cf.\ \cite[Section~12.8]{BrzWis:cor}). Note further that
$N\simeq h_C(M,C)$. Thus the notion of a ring context necessarily implies that
the right $C$-comodule $M$ is a quasi-finite injector. In the next theorem we
associate a matrix coring context to a quasi-finite injector.
\begin{theorem}\label{thm.qfi}
Let $M$ be a $(D,C)$-bicomodule and suppose that the
 right $C$-comodule $M$ is a quasi-finite injector.
Define $N:= h_C (M, C)$.  Then there exist
maps $\sigma$ and $\tau$ such that the sixtuple  $(C, D, {^C N^D}, ^D\! {M^C},
\sigma, \tau)$ is a matrix ring context.
\end{theorem}
Recall from \cite[Section~1.8]{Tak:Mor} (cf.\ \cite[Sections~12.5--12.6]{BrzWis:cor}) that if a $(D,C)$-bicomodule
is quasi-finite as a right $C$-comodule, then $h_C(M,C)$ is a
$(C,D)$-bicomodule with the left $C$-coaction
$
{}^{h_C(M,C)}\! \varrho := h_C(M,\Delta_C)$ and the right $D$-coaction
$\varrho^{h_C(M,C)}$ uniquely determined by the condition
$$
(h_C(M,C)\ot \Mro)\circ\varphi_C = (\varrho^{h_C(M,C)}\ot M)\circ \varphi_C,
$$
where $\varphi:\M^C\to h_C(M, -)\ot M$ is the unit of the adjunction.
This explains the $(D,C)$-bicomodule
structure of $N$ in the theorem. Recall further from \cite[Section~1.17]{Tak:Mor}
that for a quasi-finite right $C$-comodule $M$, the vector space $E = h_C(M,M)$
is a coalgebra with the coproduct and counit determined uniquely by relations
\begin{equation}\label{cop.coen}
(E\ot\varphi_M)\circ\varphi_M = (\Delta_E\ot M)\circ\varphi_M, \qquad (\eps_E\ot
M)\circ \varphi_M = M.
\end{equation}
$E$ is known as the {\em coendomorphism coalgebra of $M$}. Furthermore,
$M$ is an $(E,C)$-bicomodule. In addition if $M$ is a $(D,C)$-bicomodule,
then there exists a unique coalgebra map $\pi: E\to D$ such that
$\Mro= (\pi\ot M)\circ\varphi_M$ (cf.\ \cite[Section~1.18]{Tak:Mor}). Explicitly,
$\pi:= (\eps_E\ot D)\circ\varrho^E$, where $\varrho^E: E\to E\ot D$ is the right
$D$-coaction on $E$ induced by the left $D$-coaction on $M$. The strategy
for the proof of Theorem~\ref{thm.qfi} is to prove it first for $D=E$ and then
to deduce it for all $D$, using the  colagebra map $\pi: E\to D$.
\begin{lemma} \label{mmrc}
Suppose that a right $C$-comodule $M$ is a quasi-finite injector and
define $N:= h_C (M, C)$ and $E:= h_C (M, M)$.  Then there exist
maps $\sigma_E$ and $\tau_E$ such that the sixtuple $(C, E, {^C\! N^E}, ^E\! {M^C},
\sigma_E, \tau_E)$ is a matrix ring context.
\end{lemma}
\begin{proof}
First recall that $h_C(M,-)$ can be understood as a functor $\M^C\to \M^E$
which is the left adjoint to the cotensor functor $-\coten E M :\M^E\to \M^C$ (cf.\
\cite[Section~12.7]{BrzWis:cor}).
Since $M^C$ is a quasi-finite injector, $h_C(M, -) \simeq -\coten C N$ (cf.\
\cite[Section~12.8]{BrzWis:cor}). Thus there are the unit and counit
of adjunction
$$ \varphi : \M  ^C \to - \coten C N \coten E M, \qquad
 \nu : - \coten E M \coten C N \to \M ^E .$$
Define morphisms of right comodules
$$
\sigma_E : = \varphi _C : C \to C \coten C N \coten E M \simeq N \coten E M,
\qquad
 \tau_E := \nu_E : M \coten C N \simeq E \coten E M \coten C N \to E.
 $$
To see that $\sigma_E$ is a
$C$-bicomodule map use the fact that $\varphi$ is a natural tranformation
to produce commutative diagrams
$$
\xymatrix{C \ar[rr]^{\Delta_C}\ar[d]_{\varphi_C} && C\ot C\ar[d]^{\varphi_{C\ot C}} &
C \ar[rr]^{l_c}\ar[d]_{\varphi_C} && C\ot C\ar[d]^{\varphi_{C\ot C}}\\
C\coten C N\coten E M \ar[rr]_{\Delta_C\coten C N\coten E M} &&
C\ot
C\coten C N\coten E M  &
C\coten C N\coten E M \ar[rr]_{l_c\coten C N\coten E M} &&
C\ot
C\coten C N\coten E M }
$$
where $l_c(c') = c \otimes c'$, for all $c,c'\in C$.
Since $\Delta_C\coten C N \coten E M$ can be identified with the
left $C$-coaction ${}^{N \coten E M}\! \varrho$, putting these two diagrams
together we obtain, for all $c\in C$,
\begin{eqnarray*}
^{N \coten E M} \!\varrho \circ \varphi_C (c) &=& \varphi_{C \otimes C}
(c
\sw 1 \otimes c \sw 2)
= \varphi _{C \otimes C} \circ l_{c \sw 1 } (c\sw 2 )\\
&=& (l_{c \sw 1} \coten C N \coten D M) \circ \varphi_C(c \sw 2)
= c \sw 1 \otimes \varphi_C ( c \sw 2).
\end{eqnarray*}
Hence $\sigma_E = \varphi_C$ is a $C$-bicomodule map.  A similar
method can be used to show that $\tau_E$ is an $E$-bicomodule map.
By the properties of the unit and counit of adjunction,  the
composition
\begin{equation}\label{eq.iden}
C \coten C N \xrightarrow {\varphi_C \coten C N} C \coten C N \coten E M \coten C N
\xrightarrow {\nu_{C\coten C N}} C \coten C N
\end{equation}
yields the identity.  Since $\nu$ is a natural
transformation, the commutative diagrams induced by the morphisms
$\rho^N$, $l_n: E\to N\ot E$, $x\mapsto n\ot x$, and $^N \rho$, give the following equalities
\begin{eqnarray}
\nu_{N \otimes E} \circ (\varrho^N \coten E M \coten C N) = \varrho^N \circ
\nu_N \label{eq1}\\
\nu_{N \otimes E} \circ (l_n \coten E M \coten C N) = l_n \circ \nu_E
\label{eq2}
\end{eqnarray}
respectively.  Hence, for all $n \otimes m \otimes n'\in N\coten E
M\coten C N$ (summation suppressed for simplicity),
\begin{eqnarray*}
\varrho^N  \circ \nu_N (n \otimes m \otimes n')&=& \nu _{N \otimes
E} \circ (n \zz 0 \otimes n \zz 1 \otimes m \otimes n')
= \nu _{N \otimes E} \circ (l_{n\zz 0}(n \zz 1) \otimes m
\otimes n')\\
&=& l_{n \zz 0} \circ \nu_E (n \zz 1 \otimes m \otimes n')
= n \zz 0 \otimes \nu_E(n \zz 1 \otimes m \otimes n'),
\end{eqnarray*}
where the first equality is from \eqref{eq1} and last by \eqref{eq2}.
And so applying $N \otimes \eps_E$ to both sides and using the
canonical identification $N\coten E E\simeq N$, we obtain
$\nu _N  = N
\otimes \htau_E$,
where $\htau_E := \eps_E\circ\tau_E$. Since $\sigma_E = \varphi_C$, the first
of relations \eqref{eq.MRC} follows by the fact that the composition \eqref{eq.iden} is
the identity.
The other condition in \eqref{eq.MRC}
 is proven in a similar way.
\end{proof}

Note that the map $\tau_E$ constructed in the proof of Lemma~\ref{mmrc} is a
bijection, hence $(C, E, {^C\! N^E}, ^E\! {M^C},
\sigma_E, \tau_E^{-1})$ is a Morita-Takeuchi context.

In the situation of Theorem~\ref{thm.qfi}, the coalgebra map
$\pi:E\to D$ induces the map $N \coten E M\rightarrow N \coten D M$.
Using the matrix ring context $(C, E, {^C\! N^E}, ^E\! {M^C},
\sigma_E, \tau_E)$ in Lemma~\ref{mmrc}, define the required
matrix ring context $(C, D, {^C\! N^D}, ^D\! {M^C},
\sigma, \tau)$ by
$$
\sigma: C \xrightarrow{\sigma_E} N \coten E M\rightarrow N \coten D M, \qquad
\tau: M \cotc N \xrightarrow{\tau_E} E \xrightarrow{\pi} D.
$$
This completes the proof of Theorem~\ref{thm.qfi}.

The notion of a matrix ring context has a very natural
interpretation in the language of bicategories.  Consider the
bicategory of bicomodules where 0-cells are coalgebras, 1-cells
are bicomodules and 2-cells are bicomodule maps.  Define the
composite, $g \circ f$, of two 1-cells $f: X \to Y$ and $g: Y \to
Z$ to be $f \coten Y g : X \to Z$. Then there are obvious
associativity and unit isomorphisms.  When the isomorphisms
implicitly used in Definition \ref{MRC}, such as $(N \cotd M)
\cotc N \cong N \cotd (M \cotc N)$, are fully described it becomes
apparent that in this language $(C, D, g: C \to D, f:D \to C,
\sigma, \tau)$ is a matrix ring context if and only if the 2-cells
$\sigma : 1_C \Rightarrow f \circ g$ and $\tau : g \circ f
\Rightarrow 1_D$ form an adjoint pair in the bicategory.

\subsection{Infinite (firm) matrix contexts.}
The extension of comatrix coring contexts to non-unital firm rings
in \cite{GomVer:fir} (cf.\ \cite{CaeDeG:col}, both extending 
infinite comatrix corings of \cite{ElKGom:inf}) allows for a generalisation of matrix ring
contexts as in Definition~\ref{MRC} whereby one is no longer
confined to quasi-finite injectors. We outline basic properties of
such a generalisation in the present section.

Let $D$ be a non-counital coalgebra with coproduct $\Delta_D$. We say
that $D$ is a {\em firm coalgebra} if the map $\Delta_D: D\to D\coten{D}D$ is an
isomorphism. The inverse of $\Delta_D$ is denoted by $\nabla_D:D\coten{D}D
\to D$. A left (resp.\ right) non-unital comodule $M$ of a firm coalgebra $D$ is said
to be  {\em firm}, provided the coaction $\Mro: M\to D\coten{D}M$ (resp.\
$\roM : M\to M\coten{D}D$) is an isomorphism of comodules. The inverse of
coaction is denoted by ${}_M\nabla$ (resp.\ $\nabla_M$).

\begin{definition}\label{fMRC}
An {\em infinite matrix ring context}, $(C, D, {^C N^D}, ^D {M^C}, \sigma, \tau)$,
consists of a counital coalgebra $C$, firm coalgebra $D$, a $(C,D)$-bicomodule
 $N$, a $(D,C)$-bicomodule $M$, both counital as $C$-comodules and
 firm as $D$-comodules, and a pair of bicomodule maps
$$
\sigma: C \to N \coten D M, \qquad \tau: M \coten C N \to D
$$
such that the diagrams in Definition~\ref{MRC} commute.
\end{definition}
In  contrast to (finite) matrix ring context in
Definition~\ref{MRC}, the counit $\tau$ of an infinite matrix
context does not have a reduced form. Following the same line of
argument as in \cite[Theorem~1.1.3]{Cae:Bra}, one can associate a
pair of adjoint functors with any infinite matrix ring context.
\begin{proposition}\label{prop.f.adj}
Given an infinite matrix ring context $(C, D, {^C N^D}, ^D {M^C},
\sigma, \tau)$, denote by $\M^D$ the category of firm right
$D$-comodules. Then the functor $F = -\coten{C}N:\M^C\to \M^D$ is
the left adjoint of $G = -\coten{D}M :\M^D\to \M^C$.
\end{proposition}
\begin{proof}
This can be proven in the same way as Proposition~\ref{mrcaf},
provided one replaces all references to $\eps_D$ by the inverses
of the coactions such as $\nabla_Y$ etc. The unit of the
adjunction is $ \varphi: \M^C\to GF$, $\varphi_X =
(M\coten{C}\sigma)\circ\varrho^X, $ and the counit is $ \nu: FG\to
\M^C$, $\nu_Y = \nabla_Y\circ(Y\coten{D}\tau)$.
\end{proof}

Note that this adjoint pair of functors no longer extends to functors $\M^C\to \vec$,
$\vec\to \M^C$. Consequently, $M$ is no longer a quasi-finite injector as a right
$C$-comodule. Still, associated to an infinite matrix ring context are a $C$-ring and
a firm coalgebra. Their construction is very reminiscent of the construction of an
{\em elementary algebra} in the Morita theory of non-unital rings (cf.\
\cite[p.\ 36]{Cae:Bra}, \cite[p.\ 129]{Tay:big}).

\begin{proposition}\label{prop.f.ring}
Let $(C, D, {^C N^D}, {^D\! M^C}, \sigma, \tau)$ be  an infinite
matrix ring context.
\begin{zlist}
\item $\cA := N\coten{D}M$ is a $C$-ring with the product and unit
$$
\mA = (\nabla_N\coten{D}M)\circ(N\coten{D}\tau\coten{D}M) =
(N\coten{D}{}_M\!\nabla)\circ(N\coten{D}\tau\coten{D}M),
\qquad \eA = \sigma.
$$
\item $E := M\coten{C} N$ is a firm coalgebra with the coproduct
$$
\Delta_E = (M\coten{C}\sigma\coten{C}N)\circ (\roM\coten{C}N) =
(M\coten{C}\sigma\coten{C}N)\circ (M\coten{C}\Nro).
$$
\end{zlist}
\end{proposition}
\begin{proof}
(1) $\nabla_N$ is necessarily a $(C,D)$-bicomodule map, since it is the inverse of
a $(C,D)$-bicomodule map $\roN$. This means that the map $\mA$ is $C$-bicolinear.
To see that the two forms of $\mA$ are equivalent, apply $\roN\coten{D}M$ to get $N\coten{D}\tau\coten{D}M$ in both cases (note that $\roN\coten{D}M = N\coten{D}\Mro$). That $\eA$ is the unit for $\mA$ follows by commutative
diagrams in Definition~\ref{MRC}, while the associativity of $\mA$ is
clear from the definition.

(2) The map $M\coten{C}\sigma\coten{C}N$ is coassociative by the coassociativity of
coactions and colinearity of $\sigma$. Define
$$
\nabla_E :E\coten{E}E\to E, \qquad
\nabla_E = ({}_M\!\nabla \coten{C} N)\circ (\tau\coten{D}M\coten{C}N).
$$
Note that $\Delta_E$ is $D$-bicolinear,
hence, in particular $E\coten{E}E \subseteq E\coten{D}E$. Using
this, one checks that $\nabla_E$ is the inverse of $\Delta_E$ by a routine
calculation.
\end{proof}

Note that the coalgebra $E$ plays the same role as the coendomorphism
coalgebra $h_C(M,M)$ in the quasi-finite projector case.

\section{$\cA$-coendomorphism coalgebra and Galois modules}
The aim of this section is to introduce the notion of a Galois module,
to derive the structure theorem for such modules and construct the
associated Galois connection. Galois modules
are a particular class of modules of a $C$-ring $\cA$ that are quasi-finite
injectors as $C$-comodules. First we need to introduce the notion of
an $\cA$-coendomorphism coalgebra.
\subsection{The $\cA$-coendomorphism coalgebra and $C$-ring}
\begin{lemma} \label{wla}
Let $(C, D, {^C N^D}, {^D M^C}, \sigma, \tau)$ be a matrix ring
context and let $\cA$ be a $C$-ring. If $M$ is a right
$\cA$-module, via the map $\wra: M \cotc \cA \to M$, then $N$ is a
left $\cA$-module via the map
$$\wla: \cA  \coten C N \to N, \qquad \wla : = (N \otimes \htau)
\circ (N \coten D \wra \coten C N) \circ (\sigma \coten C \cA \coten C N)
\circ (^\cA\! \varrho \coten C N)$$
\end{lemma}
\begin{proof}
The map $\wla$ is left $C$-colinear because it is a composition of
left $C$-colinear maps.
We need to show that $\wla$ is
associative and unital. Throughout the proof we write
$\sigma(c) = c\su 1\ot c\su 2\in N\cotd M$ (summation assumed). The right action $\wra$ of $\cA$
on $M$ is denoted by $\lhd$ between the elements. Similarly, the map $\wla$ is
denoted by $\rhd$. In this notation
$$
\sum_i a^i\rhd n^i = \sum_i a^i\sco{-1}\su 1\htau(a^i\sco {-1}\su 2\lhd a^i\sco 0\ot n^i), \qquad \mbox{for all}\  \sum_ia^i\ot n^i\in \cA\cotc N.
$$
Take any $a\otimes a' \otimes n\in \cA\cotc \cA\cotc N$ (summation suppressed for clarity), and compute
\begin{eqnarray*}
(a\rhd(a' \rhd n))
&=&a \zz {-1} \su 1 \htau( a \zz {-1} \su 2 \lhd
a \zz 0 \otimes a' \zz {-1} \su 1  )\htau (a' \zz {-1} \su 2
\lhd a'\zz 0 \otimes n)\\
&=& a \zz {-1} \su 1 \htau( a \zz {-1} \su 2 \lhd
a \zz 0 \otimes a \zz 1 \su 1  )\htau(a \zz 1 \su 2 \lhd a'
\otimes n)\\
&=&a \zz {-1} \su 1 \htau( (a \zz {-1} \su 2 \lhd
a \zz 0) \zz 0 \otimes (a \zz {-1} \su 2 \lhd a \zz 0) \zz 1 \su
1)\\
 &&\times \htau((a \zz {-1} \su 2 \lhd a \zz 0) \zz 1 \su 2 \lhd a'
\otimes n)\\
&=& a \zz {-1} \su 1 \htau ((a \zz {-1} \su 2 \lhd a \zz 0) \lhd
a' \otimes n)
= a \zz {-1} \su 1 \htau ((a \zz {-1} \su 2 \lhd (a \zz 0a'))  \otimes n)\\
&=& (aa') \zz {-1} \su 1 \htau (((aa') \zz {-1} \su 2 \lhd
(aa')\zz 0)  \otimes n) = ((aa')\rhd n),
\end{eqnarray*}
where the second equality holds because $a\otimes a' \in \cA \cotc \cA$,
the third by the right $C$-colinearity of the right $\cA$-action on $M$.
The fourth equality comes from the second of equations
\eqref{eq.MRC}.  The fifth
equality follows because the right $\cA$-action is multiplicative
and the penultimate equality uses the colinearity of the product
$\mA : \cA \cotc \cA \to \cA$.  This
proves that the map $\wla$ is associative. The unitality of $\wla$ follows
by a similar calculation that uses the $C$-colinearity of the unit $\eA$ and
of $\sigma$, the unitality of $\wra$ and the first of  equations
\eqref{eq.MRC}.
\end{proof}

In the set-up of Lemma~\ref{wla}, the left action of matrix $C$-ring $N\coten{D}M$
on $N$ induced from  the right action described in Proposition~\ref{prop.matring}
is $N\coten{D}\htau$.
The next lemma shows that the $\cA$-actions are compatible with
the unit and counit of a matrix ring context.
\begin{lemma}\label{lemma.act}
Let $(C, D, {^C N^D}, ^D\! M^C, \sigma, \tau)$ be a matrix ring
context and let $\cA$ be a $C$-ring. Suppose that
$M$ is a right $\cA$-module, via the map $\wra: M
\cotc \cA \to M$ (denoted by $\lhd$ between elements)
and let $\wla$ be the left $\cA$-action on $N$ constructed in Lemma~\ref{wla} (denoted by $\rhd$ between elements).
\begin{zlist}
\item For all $m\ot a\ot n \in M\cotc \cA\cotc N$ (summation suppressed for
clarity),
$$
\htau(m\lhd a\otimes n) = \htau(m\ot a\rhd n).
$$
\item The following diagram
$$
\xymatrix{ \cA\ar[rr]^\roA \ar[d]_\Aro && \cA\cotc C\ar[rr]^{\cA\cotc \sigma} &&
\cA\cotc N\cotd M\ar[d]^{\wla\cotd M}\\
C\cotc \cA\ar[rr]^{\sigma\cotc \cA} && N \cotd M\cotc \cA\ar[rr]^{N\cotd \wra}&&
N\ot M}
$$
is commutative.
\end{zlist}
\end{lemma}
\begin{proof} Both statements follow by straightforward calculations  which use
the definition of $\wla$,   the second of equations
\eqref{eq.MRC} and the definition of a cotensor product (in the case of assertion (1)),
and  the $C$-colinearity of $\wra$ (in the case of assertion (2)).
\end{proof}

\begin{theorem}\label{thm.coend}
Let $\cA$ be a $C$-ring and $M$ a right $\cA$-module which is
a quasi-finite injector as a right $C$-comodule. Let $N := h_C(M,C)$, $E:= h_C(M,M)$ and define
a vector space $E_\cA(M)$ as the coequaliser
$$
\xymatrix{M \cotc \cA \cotc N
\ar@<1ex>[rr]^{\wra\cotc N}\ar@<-0.0ex>[rr]_{M\cotc \wla} && M\cotc N\ar@<.5ex>[rr]^{\pi_\cA}
&& \underset{}{E_\cA(M)} ,}
$$
where $\wra$ is the right $\cA$-action on $M$ and $\wla$ is the induced
left $\cA$-action on $N$ as in Lemma~\ref{wla} corresponding to the matrix ring context
$(C, E, {^C\! N^E}, ^E\! {M^C},
\sigma_E, \tau_E)$ in Lemma~\ref{mmrc}.
Then $E_\cA(M)$ is a coalgebra such that
$$
\xymatrix{E\simeq M\cotc N\ar@<.5ex>[rr]^{\qquad \pi_\cA}
&& \underset{}{E_\cA(M)} }
$$
is a coalgebra map. The coalgebra $E_\cA(M)$ is called an {\em $\cA$-coendomorphism coalgebra of $M$}.
\end{theorem}
\begin{proof}
Since $M^C$ is a quasi-finite injector,  $E$
is isomorphic to $M\coten C N$. The induced coproduct and counit in
$M\coten C N$ are $M\coten C \sigma_E \coten C N$ and $\htau_E$.
Lemma~\ref{lemma.act} implies that these two maps factor through
the coequaliser defining $E_\cA(M)$ and hence provide the latter with the
 coalgebra structure such that $\pi_\cA$ is a coalgebra map.
\end{proof}

\begin{corollary} \label{proof1}
Let $\cA$ be a $C$-ring and $M$ a right $\cA$-module which is
a quasi-finite injector as a $C$-comodule
and let $N:=h_C(M,C)$. Denote  the induced left $C$-coaction on
 $N$ by $\Nro$. Then
\begin{zlist}
\item $M$ is an $(E_\cA(M) , C)$ - bicomodule, with left coaction
$(\pi_\cA\ot M)\circ (M \cotc \sigma_E) \circ \roM$. Furthermore
this left coaction is right $\cA$-linear. \item $N$ is a
$(C,E_\cA(M))$ - bicomodule, with right coaction $(N\ot
\pi_\cA)\circ (\sigma_E\cotc N) \circ \Nro$. Furthermore this
right coaction is left $\cA$-linear.
\end{zlist}
\end{corollary}
\begin{proof}
That  $M$ is a bicomodule with these coactions follows immediately
from the facts that $M$ is a left comodule of $h_C(M,M)$ (with the
coaction $(M \cotc \sigma_E) \circ \roM$) and that $\pi_\cA$ in
Theorem~\ref{thm.coend} is a coalgebra map. That the left coaction
is right $\cA$-linear follows from the defining property of
$\pi_\cA$ and Lemma \ref{lemma.act}.  The second part of the
corollary is proved in a similar way.
\end{proof}

Thus to any right  $\cA$-module $M$ which is a quasi-finite injector as a
 right $C$-comodule one can associate the matrix ring context
 $(C, E_\cA(M) , {^C
N^{E_\cA(M)}}, {^{E_\cA(M)}\! M^C}, \sigma, \tau)$
 as in the proof of Theorem~\ref{thm.qfi}, i.e.\ with
$$
\sigma : C \xrightarrow{\sigma_E} N \coten E M\rightarrow N \coten {E_\cA(M)} M,
\qquad
\tau: M \cotc N \xrightarrow{\tau_E} E \xrightarrow{\pi_\cA} E_\cA(M).
$$
We refer to this context as an {\em $\cA$-coendomorphism context} associated to
$M$. The corresponding matrix $C$-ring is referred to as an
{\em $\cA$-coendomorphism
ring of $M$}. 

\subsection{Galois and principal modules}
The aim of  this subsection is to study the relationship between
$\cA$ and the $\cA$-coendomorphism
ring of $M$.

\begin{proposition}\label{prop.cocan}
Let $\cA$ be a $C$-ring and $M$ a right $\cA$-module which is a
quasi-finite injector as a $C$-comodule. Set $N:= h_C(M,C)$ and
define a map
$$
\beta : \cA \to N \otimes M, \qquad \beta := (N
\ot \wra) \circ (\sigma \cotc \cA) \circ {\Aro},
$$
where $\wra$ denotes the $\cA$-action on $M$ and
$\sigma$ is the unit of the $\cA$-coendomorphism context associated to
$M$. Write
$S$ for the coalgebra $E_\cA(M)$. Then:
\begin{zlist}
\item $\beta(\cA) \subset N \cots M$.
\item The map $\beta$ is a morphism of $C$-rings.
\end{zlist}
\end{proposition}
\begin{proof}
(1) Write $\sigma(c) = c\su 1\ot c\su 2$, for all $c\in C$.  Note
that on elements $\sigma(c) = \sigma_E(c)$, hence we use the same
notation for $\sigma_E$. Writing  $\lhd$  for the right action of
$\cA$ on $M$,  the map $\beta$ takes the following explicit form,
$\beta(a)=a \zz {-1} \su 1 \otimes a \zz {-1} \su 2 \lhd a \zz 0$,
for all $a\in \cA$.  Denote the left (resp.\ right) $S$-coaction
on $M$ (resp.\ $N$) in Corollary~\ref{proof1} by $\Mro$ (resp.\
$\roN$). Then
\begin{eqnarray*}
(N\ot \Mro)(\beta(a))
&=& a \zz {-1} \su 1 \otimes \pi_\cA (a \zz {-1} \su 2 \lhd a \zz 0 \ot a
\zz 1 \su 1) \otimes a \zz 1 \su 2\\
&=& a \zz {-1} \su 1 \otimes \pi_\cA(a \zz {-1} \su 2 \ot a \zz 0 \rhd a
\zz 1 \su 1) \otimes a \zz 1 \su 2\\
&=& a \zz {-2} \su 1 \otimes \pi_\cA(a \zz {-2} \su 2 \ot a \zz {-1} \su
1) \otimes a \zz {-1} \su 2 \lhd a \zz 0
= (\roN\ot M)(\beta(a)) ,
\end{eqnarray*}
where the first equality follows by the right $C$-colinearity of
the $\cA$-action, the second one is the defining property of $\pi_\cA$.
The third equality follows by Lemma~\ref{lemma.act}(2) and to derive
the last equality, the left $C$-colinearity of $\sigma$ was used.

(2) The map $\beta$ is left $C$-colinear by the left colinearity
of $\sigma$. It is right $C$-colinear by the right $C$-colinearity of the
$\cA$-action $\wra$.
To check that $\beta$ is a unital map, take any $c\in C$ and compute
\begin{eqnarray*}
\beta \circ \eA(c) \!\!\!\!\!&=&\!\!\!\!\! (N \cotd \wra) \circ (\sigma \cotc \cA)
\circ {^\cA \rho} \circ \eA(c)
= (N \cotd \wra) \circ (\sigma \cotc \cA)(c \sw 1 \otimes \eA(c
\sw 2))\\
&\!\!\!\!\!\!\!= &\!\!\!\!\!\!\!\!\!\! (N \cotd \wra)(c \su 1 \otimes c \su 2 \zz 0 \otimes \eA(c
\su 2 \zz 1))
\!= \!c \su 1 \otimes \wra \circ (M \cotc \eA) \circ \roM (c \su
2)
\!= \! \sigma (c),
\end{eqnarray*}
where the second equality is by the left $C$-colinearity of $\eA$,
the third equality is by the left $C$-colinearity of $\sigma$ and
the final equality is by the unitality of a right $\cA$-action.
Since $\sigma$ is the unit map for the $\cA$-coendomorphism
$C$-ring $N \cots M$, $\beta$ is a unital map as required.  A
calculation, virtually the same as that proving the associativity
of the left $\cA$-action in the proof of Lemma~\ref{wla}, confirms
that $\beta$ is a multiplicative map too. \note{It only remains to
show that $\beta$ respects the multiplicative structure. For any
$a\ot a'\in \cA\cotc \cA$ (summation suppressed),
\begin{eqnarray*}
\beta(a)\beta(a')&=& a \zz{-1} \su 1 \htau(a \zz {-1} \su 2
\lhd a \zz 0 \otimes a'\zz {-1} \su 1)\otimes a' \zz {-1} \su 2
\lhd a' \zz 0\\
&=& a \zz {-1} \su 1 \htau (a \zz {-1} \su 2 \lhd a \zz 0 \otimes
a \zz 1 \su 1) \otimes a \zz 1 \su 2 \lhd a'\\
&=& a \zz {-1} \su 1 \otimes (a \zz {-1} \su 2 \lhd a \zz 0)\lhd
a'
= a \zz {-1} \su 1 \otimes a \zz {-1} \su 2 \lhd a \zz 0 a'\\
&=& (aa') \zz {-1} \su 1 \otimes (aa')\zz {-1} \su 2 \lhd (aa')
\zz 0
= \beta (aa'),
\end{eqnarray*}
Where the second equality is because $a\otimes a' \in \cA \cotc
\cA$, the third is by the right $C$-colinearity of the action of
$\cA$ together with the first of equations \eqref{eq.MRC}
and the penultimate equality is from the left $C$-colinearity of
the multiplication map.}
\end{proof}
\begin{definition} \label{galoisdef}
Take $M \in \M _ \cA$ such that $M^C$ is a quasi-finite injector,
set $N= h_C (M, C)$ and let $S:= E_\cA(M )$ be the
$\cA$-coendomorphism coalgebra of $M$.  We say that $M$ is a {\em
Galois $\cA$-module} iff the map $\beta : \cA \to N \cots M$ in
Proposition~\ref{prop.cocan} is bijective. A Galois $\cA$-module
$M$ is said to be {\em principal} iff $M$ is
injective as a left $S$-comodule.
\end{definition}
The notion of a Galois module generalises that of a Galois $C$-ring introduced
in \cite[Section~6]{Brz:str}. To make this statement more transparent we recall a lemma and  
definition from \cite[Section~6]{Brz:str}.

\begin{lemma} \label{1-1cor}
For any $C$-ring $\cA$, there is a bijective
 correspondence between right $\cA$-actions, $\wrac
: C \cotc \cA \to C$,  and nontrivial characters $\kappa :\cA \to
k$. Here by a nontrivial character we mean a map $\kappa: \cA \to k$
which is multiplicative and satisfies $\kappa \circ \eA = \eps_C$.
\end{lemma}
\begin{proof}
The correspondence is given as follows:  given a right
$\cA$-action $\wrac$, the corresponding character is given by
$\kappa[\wrac] :=\eps_C \circ \wrac \circ\Aro$. In the
other direction, for each character $\kappa$ there is a map $\wrac
[\kappa] $ defined as $\wrac[\kappa] (c \otimes a) = \eps_C (c)
\kappa (a \zz 0)a \zz 1$.
\end{proof}

 In the case that $\cA$ has a nontrivial character, we can study the
set
$$
I_\kappa=  \{\kappa(a \zz 0) a \zz 1 - a \zz{-1} \kappa (a \zz 0)|
a \in \cA\}\subseteq C,
$$
which is easily checked to be a coideal. Hence we are able to define a
coalgebra of coinvariants $B_\kappa = C / I_\kappa$.
\begin{definition}\label{galoi.c.ring}
A  $C$-ring $\cA$ with a
nontrivial character $\kappa$ is called a {\em Galois $C$-ring} if
there exists an
isomorphism of $C$-rings $\beta: \cA \to C \coten{B_\kappa} C$ such that
$\kappa = (\eps_C \coten{B_\kappa} \eps_C) \circ \beta$.
\end{definition}

\begin{proposition} \label{cmrc}
If $C$ is a Galois module for some $C$-ring $\cA$, then $\cA$ is a
Galois $C$-ring.
\end{proposition}
\begin{proof}
Note that $C^C$ is a quasi-finite injector:
  $(C, C, {^C C^C}, ^C\! C^C, \sigma, \tau)$ is a
matrix ring context, where $\tau: C \cotc C \to C$ is the obvious isomorphism
and $\sigma = \Delta_C$,
 corresponding to the identity map $C\to C$ as in Example~\ref{ex.MRC.trivial}.
 Obviously, $C = h_C(C,C)$.   Since
$C$ has a right $\cA$-action it also has a non-trivial character
$\kappa$, provided by the 1-1 correspondence  in
Lemma~\ref{1-1cor}. In terms of this character the right $\cA$-action is, for
all $c\ot a\in C\cotc \cA$,
$\wrac(c \otimes
a)=\eps_C(c) \kappa(a \zz 0) a \zz 1$,  so, for all $a\in \cA$,
$$
\beta (a) = (N \cotd \wrac) \circ (\sigma \cotc \cA) \circ {^\cA
\rho} (a)
= a \zz {-2} \otimes a \zz {-1} \lhd a \zz 0
= a \zz {-1} \otimes \kappa( a \zz 0) a \zz 1.
$$
Hence $\kappa = (\eps_C \cotd \eps_C) \circ \beta$.  Feeding the
above explicit form of the
right $\cA$-action on $C$  into
Lemma~\ref{wla}, we obtain a left $\cA$-action on $C$,
 $\wlac (a \otimes c) = a \zz {-1} \kappa (a \zz 0) \eps_C(c)$, for all
 $a\ot c\in\cA\cotc C$.
Thus
$$
S= C \cotc C / \im (\wrac \cotc C - C \cotc \wlac)\simeq C/
\{\kappa(a \sw 0) a \sw 1 - a \sw{-1} \kappa (a \sw 0)| a \in
\cA\}=B_\kappa.
$$
Hence $\beta: \cA \to C \coten{B_\kappa} C$ makes $\cA$ into a
Galois $C$-ring.
\end{proof}

Following a similar line of argument as in \cite[Section~4.8]{Wis:gal} one proves
\begin{proposition}\label{prop.inj}
If $M$ is a right principal Galois module of a $C$-ring $\cA$, then $\cA$ is an injective
left $C$-comodule.
\end{proposition}
\begin{proof}
Suppose that $M$ is a principal Galois $\cA$-module, write $N=h_C(M,C)$ and
$S = E_\cA(M)$, and let $E=M\coten{C}N$ denote the $C$-coendomorphism
coalgebra of $M$. Since $\cA \simeq N\coten{S}M$, there is a chain
of isomorphisms
$$
N\coten{E}M\coten{C}\cA \simeq N\coten{E}M\coten{C}N\coten{S}M\simeq
N\coten{E}E\coten{S}M\simeq N\coten{S}M\simeq\cA.
$$
Explicitly the isomorphism $\cA\to N\coten{E}M\coten{C}\cA$ is
$(\sigma_E\ot\cA)\circ\Aro$, where $\sigma_E$ is the unit of the matrix
ring context in Lemma~\ref{mmrc}. Since $M$ is an injective left $S$-comodule and
$M\coten{C}\cA\simeq M\coten{C}N\coten{S}M = E\coten{S}M$, $M\coten{C}\cA$
is injective as a left $E$-comodule. Thus there exists a left $E$-comodule retraction
$p$ of the obvious inclusion $\iota: M\coten{C}\cA \to M\ot \cA$. Hence $N\coten{E}p$
is a left $C$-colinear retraction of $N\coten{E}\iota$,
 and there is a commutative diagram
with (split) exact rows
$$
\xymatrix{ 0 \ar[r] & N\coten{E}M\coten{C}\cA \ar[rr]<1ex>^{N\coten{E}\iota} && N\coten{E}M\ot \cA \ar[ll]^{N\coten{E}p}\\
0 \ar[r] & \cA \ar[u]^{\simeq} \ar[rr]^{\Aro} && C\ot\cA\ar[u]_{\sigma_E\ot\cA} \ ,
}
$$
from which a left $C$-colinear retraction of $\Aro$ is constructed.
\end{proof}

The main result of this section is contained in the following
\begin{theorem}\label{theorem1}
Let $\cA$ be a $C$-ring and $M$ a right $\cA$-module
which is a quasi-finite injector as a right $C$-comodule.  Set $N =
h_C (M, C)$  and
$S= E_\cA(M)$.  View $N \otimes M$ and $N\cots M$ as  left $\cA$-modules
with the left action as in Lemma~\ref{wla}.  Let $\beta$ be as  in
Proposition~\ref{prop.cocan}.
\begin{enumerate}
\item The following statements are equivalent
    \begin{enumerate}
    \item there exists a left $\cA$-module map $\chi: N\ot M \to \cA$ such that
    $\chi\circ\beta = \cA$ (i.e.\
    $ \beta: \cA \to N \otimes M$
    is a split monomorphism of left $\cA$-modules);
    \item M is a principal Galois $\cA$-module.
    \end{enumerate}
    \item The following statements are equivalent
\begin{enumerate}
    \item there exists a left $\cA$-module map $\hat{\chi}: N\cots M \to \cA$ such that
    $\hat{\chi}\circ\beta = \cA$ (i.e.\
    $ \beta: \cA \to N \cots M$
    is a split monomorphism of left $\cA$-modules);
    \item M is a Galois $\cA$-module.
    \end{enumerate}
\end{enumerate}
\end{theorem}
\begin{proof}

(1) (a) $\Ra$ (b) Suppose that there
exists a left $\cA$-module
retraction $\chi$ of $\beta$. This means
explicitly that, for all $a\in A$, $ \chi(\sigma (a \zz{-1}) \lhd a \zz 0)=a$, where
 $\sigma$ is the unit of the coendomorphism
ring context.    In particular,
for $a = \eA(c)$, this implies that, writing $\sigma(c) = c\su 1\ot c\su 2$,
$$
\eA(c)
= \chi (\sigma ( c \sw 1) \lhd \eA ( c \sw 2 )) \nonumber
= \chi ( c \su 1 \otimes c \su 2 \zz 0 \lhd \eA ( c \su 2 \zz
1))
= \chi \circ \sigma (c),
$$
where the first equality follows from the left $C$-colinearity of
$\eA$, the second by the right $C$-colinearity of $\sigma$ and the
third by the unitality of the $\cA$-action.  Therefore,
\begin{equation}\label{kapdel}
\chi\circ \sigma = \eA.
\end{equation}
First we prove that $M$ is an injective left $S$-module, by
constructing a left $S$-comodule retraction of the left
$S$-coaction on $M$. Define  a  map $\delta : S \otimes M \to M$
by the commutative diagram
$$
\xymatrix{ M\cotc N\ot M \ar[rr]^{\pi_\cA\ot M}\ar[d]_{M\cotc \chi} && S\ot M\ar[d]^\delta \\
M\cotc \cA \ar[rr]^\wra && M}
$$
The map $\delta$ is well defined because $\chi$ is assumed to be a
left $\cA$-module map. By equation \eqref{kapdel} and the
unitality of the right $\cA$-action we obtain, for all $m\in M$,
$$\delta \circ {\Mro(m)} = m \zz 0 \lhd \chi( \sigma (m \zz 1))
= m \zz 0 \lhd \eA(m \zz 1)= m.
$$
Hence
$\delta$ is a retraction of the left coaction. Note that $\Mro$ is
right $\cA$-linear since, for all $m\ot a \in M\cotc \cA$ (summation
suppressed),
\begin{eqnarray*}
\Mro (m \lhd a) &=&
 \pi_\cA(m \lhd a\zz 0 \ot a \zz 1 \su 1) \otimes a \zz 1 \su 2
= \pi_\cA(m \ot a \zz 0 \rhd a \zz 1 \su 1) \otimes  a \zz 1 \su
2\\
&=& \pi_\cA(m \ot a \zz {-1} \su 1) \otimes a \zz {-1} \su 2 \lhd a \zz
0
= {\Mro (m)} \lhd a,
\end{eqnarray*}
where the first equality holds because the $\cA$-action is a right
$C$-colinear.  The second equality follows by the definition of
$\pi_\cA$ and the third by Lemma~\ref{lemma.act}(2).
To derive the last equality the fact that  $m \otimes a \in M \cotc \cA$
was used.  Now it is easy to see that, for all $m\ot n\in M\cotc N$ and
$m'\in M$,
\begin{eqnarray*}
(S \otimes \delta) \circ (\Delta_S \otimes M)(\pi_\cA(m \otimes n)
\otimes m')
&=& \pi_\cA( m \zz 0 \otimes m \zz 1 \su 1) \otimes \delta (\pi_\cA(m \zz 1 \su 2
\otimes n) \otimes m')\\
&=& \pi_\cA( m \zz 0 \otimes m \zz 1 \su 1)  \otimes m \zz 1 \su 2 \lhd
\chi(n \otimes m')\\
&=& {\Mro(m)} \lhd \chi (n \otimes m')
= {\Mro} (m \lhd \chi (n \otimes m'))\\
&=& {\Mro} \circ \delta (\pi_\cA(m \otimes n) \otimes m').
\end{eqnarray*}
To understand the first equality recall that the coproduct in $S =
E_\cA(M)$ is defined as $\Delta_S(\pi_\cA(m \otimes n))= \pi_\cA(m
\zz 0 \otimes m \zz 1\su 1) \otimes \pi_\cA(m \zz 1\su 2 \otimes
n)$.  The above calculation means that $\delta$ is a left
$S$-comodule map and hence completes the proof that $M$ is an
injective $S$-comodule. Now define $\hat{\beta}=\chi | _{N \cots
M}$.  As $\im(\beta) \subset N \cots M$ and $\chi$ is a retraction
of $\beta$,
 it is clear that $\hat{\beta}
\circ \beta = \cA$.  To see that $\hat{\beta}$ is also a right inverse of
$\beta$ take an element $n \otimes m \in N \cots M$ and compute
\begin{eqnarray*}
\beta \circ \hat{\beta} (n \otimes m) &=& \sigma ( \hat{\beta}(n \otimes m)\zz
{-1}) \lhd \hat{\beta} (n \otimes m) \zz 0
= \sigma( n \zz {-1}) \lhd \hat{\beta}(n \zz 0 \otimes m)\\
&=& n \otimes m \zz 0 \lhd \hat{\beta} (\sigma (m \zz 1))
= n \otimes m \zz 0 \lhd \eA (m \zz 1)
= n \otimes m.
\end{eqnarray*}
The second equality is because $\chi$ is left
$\cA$-linear, which demands that it is left $C$-colinear.  To justify the
third equality, remember that $n \otimes m \in N\cots M$ and so,
with coactions as in Corollary~\ref{proof1},
 $(N\ot\pi_\cA\ot M)[\sigma(n \zz {-1}) \ot n \zz 0 \otimes m -n \otimes m \zz
0 \ot \sigma (m \zz 1)] = 0$.  Since $\chi$ (and hence also
$\hat{\beta}$) is left $\cA$-linear, we can apply
$(N\ot \wra)\circ (N\ot M\ot \hat{\beta})$ to this equality, thus obtaining
the third equality in the above calculation. The fourth equality follows by equation
\eqref{kapdel} and the final equality by the unitality of the right
$\cA$-action.  Thus $\hat{\beta}$ is the required inverse of $\beta$ and we
conclude that $M$ is a principal Galois $\cA$-module.

(1) (b)$\Ra$ (a)   Assume that $M$ is a principal Galois
$\cA$-module and let $\delta : S\otimes M \to M$ be an $S$-comodule
retraction of $\Mro$, i.e.,
$\delta \circ {\Mro}=M$.  We can construct a left $\cA$-linear
retraction for
$\beta$ by making the following composition
$$
\chi : N \otimes M \xrightarrow{\roN \otimes M} N \otimes
S \otimes M \xrightarrow{N \otimes \delta} N \cots M
\xrightarrow{\beta^{-1}} \cA
$$
Note that the image of the first two compositions is in $N \cots
M$ because $\delta$ is left $S$-colinear. Note further that $\chi$ is
left $\cA$-linear, since $\Nro$ is left $\cA$-linear (by an argument similar
to the proof of right $\cA$-linearity of $\roM$ in (1) (a) $\Ra$ (b)).
Furthermore
$$
\chi \circ \beta = \beta^{-1} \circ (N \otimes \delta)
\circ (\roN \otimes M) \circ \beta
= \beta^{-1} \circ (N \otimes \delta) \circ (N \otimes {\Mro}) \circ \beta
= \beta^{-1} \circ \beta
= \cA,
$$
where the second equality follows by the fact that $\im(\beta)\in N\cots M$.
Thus $\chi$ is the required retraction of $\beta$.

(2) That (b) implies (a) is obvious.
For the converse use the same method as  in the proof of the bijectivity of
$\beta$ (1) (a) $\Ra$ (b).
\end{proof}

Theorem~\ref{theorem1}, which can be understood as a dual version
of \cite[Theorem~4.4]{Brz:gal}, is the main result of the present
paper. We will we use it in Section~\ref{sec.self} to derive a
weak algebra-Galois version of Schneider's Theorem II.

\subsection{A Galois connection}
The aim of this subsection is to construct a Galois connection
associated to a matrix $C$-ring, following the method recently
employed in the case of corings in \cite{CuaGom:Gal}. Throughout
this subsection $M$ is a right $C$-comodule which is a
quasi-finite injector, $N:=h_C(M,C)$ and $E$ is the coendomorphism
coalgebra $E = h_C(M,M) \simeq M\coten{C}N$. Furthermore,  $\pi: E\to
D$ is a coalgebra epimorphism and $\cA = N\coten{D} M$ is the
associated matrix $C$-ring (cf.\ proof of Theorem~\ref {thm.qfi}).

For any subcoideal $X\subseteq \ker\pi$ (or, equivalently, a subcoextension
$E\epi E/X\epi D$) define a matrix $C$-ring
$$
\cA(X) := N\coten{E/X} M.
$$
For any subcoideal $Y\subseteq X$, the coalgebra map $E/Y\to E/X$
induces an inclusion of $C$-rings $\cA(Y)\subseteq \cA(X)$. Note
that $\cA(0) = N\coten{E}M$ and $\cA(\ker\pi) = \cA$.  In
particular $Y\subseteq \ker\pi$ induces an inclusion of $C$-rings
$\cA(Y) \subseteq \cA$.
\begin{lemma}\label{lem.inc}
For any subcoideal $X\subseteq \ker\pi$,
$$
\ker \pi_{\cA(X)} \subseteq X,
$$
where $\pi_{\cA(X)} : E\to E_{\cA(X)}(M)$ is the surjection defining the
$\cA(X)$-coendomorphism coalgebra of $M$ (cf.\ Theorem~\ref{thm.coend}).
\end{lemma}
\begin{proof}
Write $\pi_X : E\to E/X$ for the canonical coalgebra epimorphism. In view of
the form of actions of $N\coten{E/X} M$ on $M$ and $N$ in
Proposition~\ref{prop.matring} and the definition of $E_{\cA(X)}(M)$ in
Theorem~\ref{thm.coend}, $x$ is an element of $\ker \pi_{\cA(X)}$ if and
only if there exists $m\ot n\ot m'\ot n' \in M\coten{C}N\coten{E/X}M\coten{C}N$
(summation suppressed for clarity), such that
$$
x = \htau_E(m\ot n) m'\ot n' - m\ot n\htau_E(m'\ot n').
$$
Note that the $E/X$-coactions on $M$ and $N$ are
$$
\Mro(m) = \pi_X(m\zz 0\ot m\zz1\su 1)\ot m\zz 1\su 2,
\qquad \roN(n) = n\zz{-1}\su 1 \ot \pi_X(n\zz{-1}\su 2 \ot n\zz 0),
$$
where we write   $\sigma_E(c) = c\su 1\ot c\su 2$, for all $c\in C$. If
$m\ot n\ot m'\ot n' \in M\coten{C}N\coten{E/X}M\coten{C}N$, then
\begin{eqnarray*}
\htau_E(m\ot n)\pi_X(m'\ot n') &=& \htau_E(m\ot n)\pi_X(m'\zz 0\ot m'\zz1\su 1)
\htau_E(m'\zz 1\su 2\ot n') \\
&=& \htau_E(m\ot n\zz{-1}\su 1) \pi_X(n\zz{-1}\su 2 \ot n\zz 0)\htau_E(m'\ot n')\\
&=& \pi_X(m\ot n)\htau_E(m'\ot n').
\end{eqnarray*}
The first and third equalities follow by the fact that $\htau_E$
is the counit of $E$, while the second equality if a consequence
of the fact that the middle cotensor product is over $E/X$. Hence,
if $x\in \ker \pi_{\cA(X)}$, then  $x\in \ker\pi_X=X$, as required.
\end{proof}

In view of Lemma~\ref{lem.inc}, for any $C$-subring $\cB\subseteq\cA$ we
can define the subcoideal of $\ker\pi$,
$$
\cX(\cB) := \ker \pi_\cB,
$$
where $\pi_{\cB} : E\to E_{\cB}(M)$ is the surjection defining the
$\cB$-coendomorphism coalgebra of $M$ (cf.\ Theorem~\ref{thm.coend}). Note
that if $\cB\subseteq\cB'$ are $C$-subrings of $\cA$, then $\cX(\cB)\subseteq\cX(\cB')$.
Thus we have defined an order-reversing  correspondence between partially
ordered sets
$$
\xymatrix{ \{\textrm{$C$-subrings of $\cA$}\} \ar[rr]<.5ex> &&
\{\textrm{subcoideals of $\ker\pi$}\} \ar[ll]<.5ex>} ,
$$
where the subcoideals are ordered by the relation $X' \leq X$ iff
$X \subseteq X'$ and the $C$-subrings by inclusion. We now prove
that this correspondence is a Galois connection.
\begin{proposition}\label{prop.conn}
For all $C$-subrings $\cB\subseteq \cA$ and subcoideals $X\subseteq\ker\pi$,
\begin{zlist}
\item $\cB\subseteq \cA(\cX(\cB))$, and $\cB =  \cA(\cX(\cB))$ if and only if
$M$ is a Galois $\cB$-module;
\item $\cX(\cA(X)) \subseteq X$, and $\cX(\cA(X)) = X$ if and only if
$E_{\cA(X)}(M) = E/X$.
\end{zlist}
\end{proposition}
\begin{proof} (1) Compute,
$$
\cA(\cX(\cB)) = \cA(\ker\pi_\cB) = N\coten{E/\ker\pi_\cB} M = N \coten{E_\cB(M)}M.
$$
By Lemma~\ref{lem.inc}, there is a coalgebra map $E_\cB(M)\to D$, and we can
consider the following commutative diagram with exact rows.
$$
\xymatrix{ 0 \ar[r]  & N \coten{E_\cB(M)}M \ar[rr] && N\coten{D} M \ar@{=}[d] \\
0 \ar[r] & \cB \ar[u]^\beta \ar[rr] && \cA},
$$
where $\beta$ is the map in Proposition~\ref{prop.cocan} (with $\cB$ in place of
$\cA$).  An easy calculation reveals that,
for all $b\in \cB$, $\beta(b) =b$. Therefore, the diagram is commutative and
 $\beta$ is  the required inclusion.
 By the definition of a Galois $\cB$-module, the map $\beta$ is  identity iff
 $M$ is Galois.

 (2) Note that $\cX(\cA(X)) = \ker \pi_{\cA(X)}$ and the assertion follows by
 Lemma~\ref{lem.inc}
\end{proof}
\begin{remark}\label{rem.misgal}
By setting $\cA=\cB$ in the diagram in the proof of the first
part of Proposition~\ref{prop.conn}, it is immediately apparent
that $M$ is a Galois $\cA$-module.  Moreover this is true for any
matrix $C$-ring arising naturally from a coalgebra epimorphism
with domain $E$ (cf. proof of Theorem \ref{thm.qfi}).
\end{remark}
\begin{corollary}
The Galois connection constructed in Proposition~\ref{prop.conn}
establishes a one-to-one correspondence between $C$-subrings $\cB
\subseteq \cA$ such that $M$ is Galois $\cB$-module and
subcoideals $X \subseteq \ker\pi$ such that $E_{\cA(X)} (M) =
E/X$.
\end{corollary}
\begin{proof}
For any subcoideal $X\subseteq \ker\pi$, $M$ is a Galois
$\cA(X)$-module by Remark~\ref{rem.misgal}.  On the other hand if
$M$ is a Galois $\cB$-module, then $\cX(\cA(\cX(\cB)))=\cX(\cB)$ by
the first part of Proposition~\ref{prop.conn}. Therefore
$\cX(\cB)$ is a subcoideal of $\ker\pi$ satisfying the required
property by the second part of Proposition~\ref{prop.conn}.
\end{proof}

The Galois connection constructed in Proposition~\ref{prop.conn}
establishes a correspondence between `intermediate coextensions'
$E\epi B\epi D$ and sub $C$-rings $\cB\subseteq \cA$ and can be
understood as a dual version of the Galois connection for comatrix
corings described in \cite[Proposition~2.1]{CuaGom:Gal}. The
latter is a generalisation of a Galois connection for Sweedler
corings introduced in \cite[Proposition~6.1]{Sch:big} as a
straightforward extension of the correspondence in Sweedler's
Fundamental Theorem \cite[Theorem~2.1]{Swe:pre}.

\section{$C$-rings associated to invertible weak entwining structures.}

 Recall from \cite{CaeDeG:mod} that a {\em (right-right) weak entwining structure} is a
triple $(A,C,\psi_R)$, where $A$ is an algebra, $C$ a coalgebra,
and $\psi_R : C \otimes A \to A \otimes C$ a k-linear map which,
writing, $\psi_R(c\otimes a)=\suma a \sba \otimes c \spa$,
$(A\ot\psi_R)\circ (\psi_R\ot A)(c\ot a\ot b) = \sumab a\sba\ot
b\sbb \ot c\spab$, etc., satisfies the relations
    \begin{eqnarray}
    \suma (ab)\sba \otimes c\spa = \sumab a \sba b \sbb \otimes
    c\spab , \label{WE1}\\
    \suma a \sba \cuc(c\spa)=\suma \cuc(c\spa)1\sba a , \label{WE2}\\
    \suma a \sba \otimes \cmc (c \spa) = \sumab a \sbab \otimes
    c\sw 1 \spb \otimes c\sw 2\spa , \label{WE3}\\
    \suma 1 \sba \otimes c \spa = \suma \cuc(c \sw 1 \spa)1\sba
    \otimes c\sw 2. \label{WE4}
    \end{eqnarray}
This is a generalisation of the notion of a (right-right) entwining structure \cite{BrzMaj:coa}, motivated by the representation theory of weak Hopf algebras (cf.\ \cite{Boh:wea}, \cite{Boh:Doi}). Associated to a weak entwining structure $(A,C,\psi_R)$ is the category
$\rema$ of {\em right weak entwined modules}, i.e.\
 vector spaces $M$ together with a right $A$-action $\varrho_M$ and a
    right $C$-coaction $\roM$ such that
    \begin{equation}
    \roM \circ \varrho_M= (\varrho_M\ot C)\circ
    (M\ot\re)\circ(\roM \otimes A)
    \label{rem}
    \end{equation}
 Also associated to a (right-right) entwining structure $(A,C,\psi_R)$ are projections
 \begin{eqnarray}
 \lpr: C \otimes A \to C \otimes A, \qquad \lpr = (C \otimes
A \otimes \eps_C) \circ (C \otimes \re) \circ (\Delta_C \otimes A) , \label{prbar}  \\
p_R : A \otimes C \to A \otimes C, \qquad p_R=(\mu_A \otimes C)
\circ (A \otimes \re) \circ (A \otimes C \otimes 1_A). \label{pr}
\end{eqnarray}
That these are projections follows by equations \eqref{WE3} (in the case of $\lpr$) and
\eqref{WE1} (in the case of $p_R$). Note further that
\begin{equation}\label{projr}
\re\circ\lpr = p_R\circ\re = \re.
\end{equation}
As explained in \cite{Brz:str}, the projection $p_R$ can be used
to associate an $A$-coring to a weak entwining structure. On the
other hand, $\lpr$ is needed for associating a $C$-ring to $(A,C,
\psi_R)$ as follows:

\begin{theorem}\label{thm.weak.C-ring}
Let $(A,C,\psi_R)$ be a (right-right) weak entwining structure and let
$$ \cA = \im\ \lpr = \{ \sum_{\alpha, i} c^i\sw 1 \tens a^i_\alpha \eps_C(c^i \sw 2
^\alpha)\; |\; \sum_i a^i \tens c^i \in A \tens C\}.
$$
Then:
\begin{zlist}
\item $\cA$ is a $(C,C)$-bicomodule with the left coaction $\Aro :=
\Delta_C \tens A$ and the right coaction $\rho^{\cA}:=(C \tens \psi_R)
\circ (\Delta_C \tens A)$.
\item  The $(C,C)$-bicomodule $\cA$  is a $C$-ring with product
$$\mA : \cA \cotc \cA \to \cA, \qquad \sum_i c_i \otimes a_i \otimes c'_i \otimes
a'_i \mapsto \sum_i c_i \otimes \eps_C(c'_i) a_ia'_i,
$$
and unit
$$
\eA:C \to \cA, \qquad c\mapsto \lpr (c\otimes 1).
$$
\item
$\M _\cA \equiv \rema$.
\end{zlist}
\end{theorem}
\begin{proof}
(1)
That $\Aro$ is a left coaction follows immediately from
properties of the comultiplication.  For $\roA$,  using \eqref{WE3} note
that, for all $a\in A$ and $c\in C$,
\begin{equation}
\suma \roA (c \sw 1 \otimes a \sba \eps_C (c \sw 2 \spa)) 
=\suma c \sw 1 \otimes a \sba \otimes c \sw 2 \spa . \label{*}
\end{equation}
We aim to show
that $\roA\circ\lpr= (\lpr \otimes C) \circ \roA \circ\lpr $; because
$\lpr$ is a projection, this will imply that $\roA(\cA)\subset
\cA \otimes C$. Applying $\lpr \otimes C$ to (\ref{*}) we obtain
\begin{eqnarray*}
\suma (\lpr \otimes C) \circ \roA \circ \lpr(c\ot a)
&=& \sumab c \sw 1 \otimes a \sbab \eps_C (c \sw 2 \spb) \otimes
c\sw 3 \spa \\
&=& \suma c \sw 1 \otimes a \sba \otimes c \sw 2 \spa
= \rho ^ \cA \circ \lpr(c\ot a),
\end{eqnarray*}
where the second equality is by (\ref{WE3}) and the third  by
(\ref{*}).  To see that $\roA$ is counital simply apply
$C\otimes A \otimes \eps_C$ to (\ref{*}).  To complete the proof that
$\roA$ is a coaction only remains to prove that it is
coassociative.  Take any $c\ot a\in C\ot A$ and compute
\begin{eqnarray*}
(\rho^\cA \otimes C) \circ \rho^\cA\circ \lpr(c\ot a)
\!\!\!&=&\!\!\! \suma \roA (c \sw 1 \otimes a \sba ) \otimes c
\sw 2 \spa
= \sumab c\sw 1 \otimes a \sbab \otimes c \sw 2 \spb \otimes c \sw 3
\spa\\
&&\hspace{-1in} = \suma c \sw 1 \otimes a \sba \otimes c \sw 2 \spa \sw 1 \otimes c
\sw 2 \spa \sw 2
= (\cA \otimes \Delta_C) \circ \rho^\cA\circ \lpr(c\ot a),
\end{eqnarray*}
where the first and last  equalities follow by (\ref{*}) and the
third by (\ref{WE3}).  Using the coassociativity of the coproduct one easily checks that
left and right coactions commute with each other, thus making $\cA$ into a
$(C,C)$-bicomodule, as claimed.

(2)
The map $\mA$ is obviously left $C$-colinear. A simple  calculation,
which uses \eqref{*}, confirms that $\mA$ is also a right $C$-comodule
map. Similarly, $\eA$ is obviously left $C$-colinear. Using \eqref{WE4} and
\eqref{*} we immediately find
$$
\lpr (c\sw 1\ot 1) \ot c\sw 2 = \suma c\sw 1\ot 1\sba\ot c\sw 2\spa = \roA\circ\lpr(1\ot c),
$$
hence $\eA$ is right $C$-colinear as well. A straightforward
calculation proves that $\mA$ is associative and unital.

(3) Let $\Psi: \M _\cA \to \rema$ be the map which leaves each
$\cA$-module unchanged as a  $C$-comodule, but which changes
the right $\cA$-action $\wra$ into a map $\Psi(\wra):M \otimes A \to M$,
 $\Psi(\wra) = \wra \circ ( M \cotc \lpr) \circ (\roM \otimes
A)$, which will presently be shown to be a right $A$-action for which $M$
is an entwined module. Unitality follows easily as
$$
\Psi(\wra)(m\otimes 1) = \wra \circ (M \cotc \lpr) (m \zz 0 \otimes
m \zz 1 \otimes 1)
= \wra \circ (M \cotc \eA) \circ \roM (m)
= m,
$$
where the second equality is by the definition of the unit and
last equality comes from the unitality of an $\cA$-action. For associativity,
take any $a,a'\in A$ and $m\in M$ and compute
\begin{eqnarray*}
\Psi(\wra) &&\hspace{-.4in} \circ (\Psi(\wra) \otimes A) (m \otimes a \otimes a')\\
&=& \suma \wra \circ (M \cotc \lpr) (\roM \circ \wra (m \zz 0 \ot m \zz
1
\otimes a \sba \eps_C (m \zz 2 \spa)) \otimes a')\\
&=&\suma  \wra \circ (M \cotc \lpr) ((\wra \otimes C)(m \zz 0 \ot m \zz 1
 \otimes a \sba  \otimes m \zz 2 \spa) \otimes a')\\
&=& \suma \wra \circ (\wra \cotc \mathcal{A})(m\zz 0 \ot (m \zz 1
\otimes a \sba) \ot \lpr (m \zz 2 \spa \otimes a'))\\
&=& \suma \wra \circ (M \cotc \mA)(m\zz 0 \ot (m \zz 1 \otimes a \sba)
\ot \lpr (m \zz 2 \spa \otimes a'))\\
&=& \sumab \wra(m \zz 0 \ot m \zz 1 \otimes a \sba a'\sbb \eps_C(m \zz 2
\spa \sw 1) \eps_C(m\zz 2 \spa \sw2 \spb))\\
&=& \suma \wra ( m\zz 0 \ot m \zz 1 \otimes (aa')\sba \eps_C (m \zz 2
\spa))
= \Psi(\wra)(m \otimes a  a').
\end{eqnarray*}
Here the second equality is from the right $C$-colinearity of the
map $\wra$ and the equality (\ref{*}).  The fourth equality comes
from the associativity of $\wra$.
The penultimate equality follows from the
 definition of a counit and  (\ref{WE1}).
Next we check that this right action makes $M$ an entwined module:
\begin{eqnarray*}
( \Psi(\wra)\ot C)\circ (M\ot\psi_R)\circ(\roM\ot A) &=& (\wra\ot C)\circ (M\ot C\ot \psi_R)\circ (M\ot\Delta_C\ot A)\circ (\roM\ot A)\\
&=&  \wra \circ (M \cotc \roA\circ \lpr) \circ (\roM \otimes A)\\
&=& \roM\circ \wra \circ ( M \cotc \lpr) \circ (\roM \otimes
A)
= \roM \circ \Psi(\wra),
\end{eqnarray*}
where the first  equality follows by the coassociativity of a coaction,
the definition of a counit and (\ref{WE3}),
the second by
(\ref{*}) and penultimate equality by the colinearity of $\wra$.

Given a morphism $f: M\to N$ in $\M_\cA$, we define $\Psi(f) = f$. Using the
$C$-colinearity of $f$ and that $\overline{\rho_N}\circ (f \cotc \cA)= f\circ \wra$,
one easily finds that the map $f$ is also right $A$-linear, when $M$ and $N$ are
viewed as $A$-modules with actions $\Psi(\wra)$ and $\Psi(\overline{\rho_N})$
respectively. Thus $\Psi$ is a functor.

In the other direction, define $\Theta:\rem \to \M _\cA$ to be the
map which leaves each entwined module $M$ unchanged as a
$C$-comodule, but which changes the right $A$-action $\rho_M$
into a map $\Theta(\rho_M):M \cotc \cA \to M$ defined as
$\Theta(\rho_M)= \rho_M \circ (M \otimes \eps_C \otimes A)$. Since
$M\cotc \cA = (M\cotc \lpr)(M\cotc C\ot A)$, all elements of
$M\cotc \cA$ are linear combinations of $x = \suma m\zz 0 \ot m\zz
1\ot a\sba\eps_C(m\zz 2\spa)$ with $a\in A$ and $m\in M$. In  view
of the fact that $M$ is an entwined module, $\Theta(\rho_M) (x)=
ma$. From this, the unitality and associativity of
$\Theta(\rho_M)$ easily follow. The right $C$-colinearity of
$\Theta(\rho_M)$ is confirmed by the following simple calculation
that uses that $M$ is an entwined module and equation \eqref{*}:
\begin{eqnarray*}
\roM \circ \Theta ( \rho_M )(x) &=& \suma m\zz 0a\sba \ot m\zz 1\spa
=
(\rho_M \ot C)\circ (M\ot \eps_C\ot A\ot C)
\circ (M\ot \roA)(x)\\
&= &(\Theta ( \rho_M )\ot C)\circ (M\ot \roA)(x).
\end{eqnarray*}

Given a morphism $f: M\to N$ in $\rem$, define $\Theta(f) = f$.
Then $\Theta(f)$ is obviously right $C$-colinear and is right
$\cA$-linear by the definition of the $\cA$-action and the
$A$-linearity of $f$. Since the composition in both categories is
provided by the composition in the category of vector spaces,
$\Theta: \rem\to \M_\cA$ is a functor.

Finally, note that for all $M\in\rem$, $m\in M$ and $a\in A$,
$$
\Psi(\Theta(\rho_M))(m\ot a) = \suma \rho_M(m\zz 0\ot a\sba \eps_C(m\zz 1^\alpha))
= (ma)\zz 0\eps_C((ma)\zz 1) = \rho_M(m\ot a),
$$
where the second equality follows by the fact that $M$ is an entwined module. On the
other hand, taking $M\in \M_\cA$ and applying $(\Theta(\Psi(\wra))$ to
$x = \suma m\zz 0 \ot m\zz 1\ot a\sba\eps_C(m\zz 2\spa)$ one immediately obtains
that $(\Theta(\Psi(\wra))(x)  = \wra(x)$. Therefore, $\Psi$ and $\Theta$ are inverse
isomorphisms of the categories, as required.
\end{proof}

As explained in \cite[Example~2.4]{Brz:str}, there is a weak
entwining structure associated to any weak coalgebra-Galois
extension. Dually, there is a weak entwining structure associated
to a {\em weak algebra-Galois coextension} as described in the
following
\begin{example}\label{crwce}
Let $A$ be an algebra, $C$ be a coalgebra and a right $A$-module
with the action $ \roC$.  Define the coideal
$$
I= \{(ca) \sw 1 \alpha((ca)\sw 2)- c \sw 1 \alpha (c \sw 2 a)
| a \in A, c \in C, \alpha \in \textrm{Hom}(C , k)\},
$$
let $B=C /I$ and let
$$
\overline{\beta}: C \otimes A \to C \cotb C, \qquad
\overline{\beta}:= (C \otimes {\roC} ) \circ (\Delta_C \otimes
A).
$$
View $C \cotb C$ as an object of $^C \M _A$ in the obvious way.
Now suppose that $C \twoheadrightarrow B$ is a {\em weak
algebra-Galois coextension}, i.e.\ that there exists a morphism
$\overline{\chi}: C \cotb C \to C \otimes A$ in $^C \M _A$ such
that $\overline{\beta} \circ \overline{\chi} = C \cotb C$. Let
$\omega : C \cotb C \to A$, $\omega := (\eps_C \otimes A) \circ
\overline{\chi}$ be the cotranslation map. Define
$$
\re : C \otimes A \to A \otimes C, \qquad \re:= (\omega \otimes C)
\circ (C \otimes \Delta_C) \circ \overline{\beta}.
$$
Then $(A, C, \re)$ is a (right-right) weak entwining structure.
Moreover $\re$ is the unique weak entwining map such that $C \in
\rema$ with structure maps $\Delta_C$ and $\roC$. This example can
be proven along the same lines as  \cite[Theorem~3.5]{BrzHaj:coa}.
\end{example}
\begin{remark}\label{BequalBk}
If $C \in \rema$, then the definition of $I$ coincides with that
of $I_\kappa$ in the definition of a Galois $C$-ring
(Definition~\ref{galoi.c.ring}), where $\kappa$ is the restriction
of $\eps_C\circ\roC$ to $\cA$.
\end{remark}
\begin{remark}\label{entg}
If $\cA$ is a $C$-ring associated to a weak entwining structure,
then $\cA$ is a Galois $C$-ring iff $\overline{\beta}|_\cA : \cA
\to C \cotb C$ is a bijection.
\end{remark}
A connection between weak algebra-Galois coextensions and Galois
$C$-rings (hence also Galois $\cA$-modules) is provided by the
following
\begin{proposition}\label{prop.Gal}
The $C$-ring associated to the weak entwining structure in
Example~\ref{crwce} is a Galois $C$-ring. Conversely, if the
$C$-ring associated to a weak entwining structure $(A,C,\psi_R)$
is a Galois $C$-ring, then $C$ is a weak algebra-Galois
coextension.
\end{proposition}
\begin{proof}
If $\cA$ is the $C$-ring associated to the weak
entwining structure in Example~\ref{crwce}, then $\cA= \mathrm{Im}
(\overline{\chi} \circ \overline{\beta})$.  Since
$\overline{\beta} \circ \overline{\chi} = C \cotb
C$, the map $\overline{\beta}\mid_\cA$ is a bijection.
Therefore, by Remark~\ref{entg}, $\cA$ is a Galois $C$-ring.
Conversely if $\cA$ is a Galois $C$-ring  and associated to a weak
entwining structure then by Remark~\ref{entg},
$\overline{\beta}|_\cA : \cA \to C \cotb C$ is a bijection,
furthermore it is clear from the definition of $\overline{\beta}$
that it is a morphism in $^C \M _A$.  Now observe that the
composition of the maps ${\overline{\beta}|_\cA}^{-1}:C \cotb C
\to \cA$ and then the inclusion $\cA \hookrightarrow C \otimes A$
is a morphism in $^C \M _A$ splitting $\overline{\beta}$.
Therefore $C \twoheadrightarrow B$ is a weak algebra-Galois
coextension.
\end{proof}

The notion of a  (right-right) weak entwining structure has a
left-handed counterpart. A {\em (left-left) weak entwining
structure} is a triple
    $(A,C,\le)$ consisting of an algebra $A$, a coalgebra $C$, and a $k$-linear map
    $\le : A \otimes C \to C \otimes A$ which, writing,
    $\le(a\otimes c)=\sumA c _{E} \otimes a ^{E} $, $\le(a\otimes c)=\sum_F c _{F} \otimes a ^{F} $ etc., satisfies the relations
    \begin{eqnarray}
    \sumA c_{E} \otimes (ab)^{E}=\sumAB c_{EF}\otimes a^{F}b^{E},  \label{le1}\\
    \sumA \cuc(c_{E})a^{E} = \sumA a \cuc(c_{E})1^{E},  \label{le2}\\
    \sumA \cmc(c_{E}) \otimes a^{E} = \sumAB c _{(1)E} \otimes c_{(2)F} \otimes
a^{EF},  \label{le3}\\
    \sumA c_{E} \otimes 1^{E}=\sumA c _{(1)}\otimes \cuc
    (c_{(2)E})1^{E}.  \label{le4}
    \end{eqnarray}
    Associated to a (left-left) entwining structure is the category
    of left entwined modules $\lem$ defined by the obvious modification of
    condition \eqref{rem}. Also, there are projections
\begin{eqnarray}
\lpl: A \otimes C \to A \otimes C, \qquad \lpl = (\eps_C \otimes A
\otimes C) \circ (\le \otimes C) \circ (A \otimes \Delta_C), \label{plbar}\\
p_L : C \otimes A \to C \otimes A, \qquad p_L=(C \otimes \mu_A)
\circ (\le \otimes A) \circ ( 1 \otimes C \otimes A). \label{pl}
\end{eqnarray}
Note that
\begin{equation}\label{projl}
\le\circ\lpl = p_L\circ\le = \le.
\end{equation}
In an analogous  way as in Theorem~\ref{thm.weak.C-ring}, $\cB =
\im\; \lpl$ is a $C$-ring, and $\lem\equiv {}_\cB\M$.  Note that
the left and right $C$-coactions on $\cB$ are given by $\Bro =
(\le\ot C)\circ (A\ot\Delta_C)$, $\roB = A\ot\Delta_C$,
respectively. In the case of invertible weak entwining structures
the $C$-rings associated to the left and right weak entwining
structures are strictly related. Recall from \cite{BrzTur:str}
\begin{definition}\label{def.inv}
 An {\em invertible weak entwining structure} is a quadruple
$(A,C, \re, \le)$ such that
\begin{blist}
 \item $(A,C,\re)$ is a right-right weak entwining structure and $(A,C,\le)$ is a
 left-left weak entwining structure;
 \item $
\re \circ \le = p_R$ and $ \le \circ \re = p_L$.
\end{blist}
\end{definition}
As observed in \cite{AloFer:inv}, if $(A,C, \re, \le)$ is an invertible weak entwining
structure, then for all $c\in C$,
\begin{equation}\label{cond.c}
\sumA\eps_C(c_E)1^E=\suma 1\sba \eps_C (c\spa
).
 \end{equation}
\begin{lemma}[cf.\  Proposition~1.5 in \cite{AloFer:inv}]\label{lemma1}
Let $(A, C, \re, \le)$ be an invertible weak entwining structure.
Then
$\lpr = p_L$ and $\lpl = p_R$.
\end{lemma}
\begin{proof}
To see that $\lpr=p_L$, take any $a\in A$ and $c\in C$, and compute
\begin{eqnarray*}
\lpr(c \otimes a) &=& \suma c \sw 1\otimes a \sba \eps_C (c \sw 2 \spa)
= \suma c \sw 1 \otimes \eps_C (c \sw 2 \spa) 1 \sba a\\
&=& \sumA c \sw 1 \otimes \eps_C (c \sw 2 _E)1^E a
= \sumA c_E \otimes 1^E a
= p_L (c \otimes a ),
\end{eqnarray*}
where the second equality follows by \eqref{WE2}, the third by \eqref{cond.c},
 and the fourth by \eqref{le4}.  A similar
calculation shows that $\lpl = p_R$.
\end{proof}
\begin{remark} \label{bstar}
 Lemma~\ref{lemma1} shows
that conditions (b) in the definition of an
invertible weak entwining structure may be replaced with
alternative conditions:\medskip

 (b*)  $\re \circ \le =
\lpl$ and $ \le \circ \re = \lpr$.\medskip

\noindent Note further that both $\cA$ and $\cB$ are not only $C$-rings but
also $A$-corings.
\end{remark}
\begin{proposition}\label{isom}
Let $(A, C, \re, \le)$ be an invertible weak entwining structure
and let $\cA = \mathrm{Im }\ \lpr$ and $\cB = \mathrm{Im }\ \lpl$
be the corresponding $C$-rings.  Then the restrictions of the
entwining maps
$$ \le : \cB \to \cA, \qquad \re : \cA \to \cB $$
are inverse isomorphisms of $C$-rings.
\end{proposition}
\begin{proof}
Since $\lpr$ and $\lpl$ are projections, the conditions (b*) in
Remark~\ref{bstar} imply that the restrictions of $\re$ and $\le$
to $\im \ \lpr$ and $\im \ \lpl$ respectively, are inverse
isomorphisms of vector spaces.  Using \eqref{WE3} one easily finds that
$$
(\re \otimes C) \circ \roA\circ \lpr = (A\ot \Delta_C)\circ \re\circ\lpr = \roB \circ\re\circ\lpr,
$$
where the second equality follows by the definition
of the right $C$-coaction on $\cB$.
This shows
that $\re$ is right $C$-colinear.  Similarly to show the $\re$ is
left $C$-colinear compute
\begin{eqnarray*}
\Bro\circ\re\circ\lpr &=& (\le\ot C)\circ(\re\ot C)\circ(C\ot\re)\circ
(\Delta_C\ot A)\circ\lpr\\
&=& (\lpr \ot C)\circ\roA\circ\lpr = \roA\circ\lpr = (C\ot \re)\circ \Aro\circ\lpr,
\end{eqnarray*}
where the first equality follows by the definition of $\Bro$ and
property \eqref{WE3} and the second by the definition of an
invertible weak entwining structure and the definition of the
coaction $\Aro$. The third is a consequence of the fact that the
image of $\rho^\cA$ is in $\cA\ot C$ (compare the proof of
Theorem~\ref{thm.weak.C-ring}(1)), and the last equality is
immediate from the definitions of $\Aro$ and $\roA$ in
Theorem~\ref{thm.weak.C-ring}(1).
Hence
$\re$ is a $(C,C)$-bicomodule map.  Similarly one shows that $\le$
is a $(C, C)$-bicomodule map.  The unitality of $\re$ is easily
checked with the help of Lemma~\ref{lemma1},  \eqref{WE4} and
\eqref{cond.c},
$$
\re \circ \eA (c) = \re \circ \lpr (c\otimes 1)
= \re(c\ot 1)= \suma \eps_C ( c \sw 1\spa) 1 \sba \otimes c \sw 2
= \lpl (1 \otimes c)
= \eB(c).
$$ Since $\cA\cotc \cA = (\lpr\cotc \lpr)(C\ot A\cotc C\ot A)$, it suffices
to check the multiplicativity of $\re$ on elements of the form
$$x = \suma\lpr(c\sw 1\ot a\sba)\ot\lpr(c\sw 2\spa\ot a') =
 \sumab c\sw 1\ot a\sbb\ot c\sw 2\spb\sw 1\ot a'\sba\eps_C(c\sw 2\spb
 \sw 2\spa).
 $$
 The definition of product in $\cA$ and properties \eqref{WE1} and \eqref{projr}
 yield
 $$
 \re \circ \mA (x) = \sumab \re (c\sw 1 \ot a\sbb a'\sba\eps_C(c\sw 2^{\beta\alpha}))
 =\re\circ\lpr (c\ot aa') = \re(c\ot aa').
 $$
 On the other hand, in view of \eqref{projr} and conditions
 \eqref{WE1} and \eqref{WE3}
 \begin{eqnarray*}
 \mB\circ(\re\cotc\re)(x) &= & \suma \mB\circ(\re\cotc\re)(c\sw 1\ot a\sba\ot c\sw 2\spa\ot a') \\
 &=& \sum_{\alpha,\beta,\gamma}a_{\alpha\beta}a'_\gamma\eps_C(c\sw 1\spb)\ot
 c\sw 2^{\alpha\gamma} = \re(c\ot aa').
 \end{eqnarray*}
 Thus $\re$ is multiplicative, hence a $C$-ring isomorphism as required.
\end{proof}
\begin{corollary} \label{cla}
Let $(A, C, \re, \le)$ be an invertible weak entwining structure.
If $C \in \rema$, then $C \in \lem$ with the action, for all $a\in A$, $c\in C$,
$$
a c = \sumA  c_E \sw 1\eps_C(c _E \sw 2  a ^E) .
$$
\end{corollary}
\begin{proof}
To see this make the following chain of deductions.  First, if $C
\in \rema$, then $C \in \M _\cA$ by Theorem~\ref{thm.weak.C-ring}.
The corresponding right $\cA$-action is, for all $c\otimes a\in
\cA$ (summation suppressed
 for clarity) and $c'\in C$,
$$
c' \lhd ( c \otimes a) = \eps_C (c) c' a.
 $$
Since there is an obvious matrix ring context $(C, C, {^C
C^C}, ^C {C^C}, \sigma, \tau)$ (cf.\ Example~\ref{ex.MRC.trivial} or the proof of Proposition~\ref{cmrc}), by Lemma~\ref{wla}  $C$ is a left $\cA$-module   with left $\cA$-action
$$ (c \otimes a) \rhd c' = c \sw 1 \eps_C (c \sw 2  a) \eps_C
(c').
$$
By Proposition~\ref{isom}, $\le : \cB \to \cA$ is an
isomorphism of $C$-rings and so $C \in { _\cB \M }$ with
left $\cB$-action
$$
( a \otimes c) \rhd c' = \sumA ( c_E \otimes a ^E) \rhd c' = \sumA c _E \sw 1
\eps_C (c_E \sw 2 a ^E) \eps_C (c').
$$
Finally we use the
correspondence $_\cB \M \equiv \lem$ to view $C$ in $\lem$ with
the left $A$-action as stated.
\end{proof}

\section{Coextensions of self-injective algebras}\label{sec.self}
In this section we start with an invertible weak entwining
structure such that $C$ is a right entwined module and then use
 Theorem~\ref{theorem1} to deduce a
criterion for this coalgebra to be a weak $A$-Galois coextension.
Since we will work in this setting, $S= E_\cA (C)$ (where $\cA$ is
the $C$-ring associated to the (right-right) weak entwining
structure) will be the same as $B_\kappa$, by the isomorphism
given in the proof of Proposition~\ref{cmrc}.  Moreover, as stated
in Remark~\ref{BequalBk}, $B_\kappa=B$ so for simplicity we shall
henceforth denote all these objects by $B$.
\begin{proposition}\label{prop.co.ext}
Let $(A , C, \re, \le)$ be an invertible weak entwining structure
such that $C$ is a right entwined module,
and let $\cA$ be the $C$-ring corresponding to $(A, C, \re)$.
View $C$ as a left $A$-module as in Corollary~\ref{cla}. Then
$ C \twoheadrightarrow B $ is a weak $A$-Galois coextension and
$C$ is  injective as a left $B$-comodule if and only if there exists
a $k$-linear map $\hat{g}:C \otimes C \to A$ such that, for all
$c\in C$ and $a\in A$,
\begin{equation} \label{cond}
\sum_\alpha a \sba \hat{g} (c \spa \otimes c')= \sum \sba
\hat{g}(a \sba c \spa \otimes c'),
\end{equation}
and
\begin{equation} \label{cotran}
\hat{g}(c\sw 1\ot c\sw 2a) = \suma a\sba\eps_C(c\spa).
\end{equation}
\end{proposition}
Since it is assumed in Proposition~\ref{prop.co.ext} that $C$ is a
weak entwined module with an $A$-action $\roC$, $C$ is a right
$\cA$-module. In view of Proposition~\ref{prop.Gal} and
Proposition~\ref{cmrc}, to prove Proposition~\ref{prop.co.ext} we
need to find criteria for $C$ to be a principal Galois $\cA$-module.
Setting $\kappa = \eps_C\circ\roC$ in the construction of the
proof of Proposition~\ref{cmrc}, we obtain  a left $\cA$-module
structure on $C\ot C$,
$$ \overline{_{C\otimes C}\rho} : \cA \otimes C\simeq \cA\coten{C}C\ot C
 \to C \otimes C, \qquad c
\otimes a \otimes c' \mapsto c \sw 1 \eps_C (c \sw 2  a)
\otimes c'.$$
In view of  Theorem~\ref{theorem1} we need to study $\cA$-module
retractions of $\beta$ (or $\overline{\beta}$). First we classify all candidates for such retractions.

\begin{lemma}\label{lemma2}
Given an invertible weak entwining structure $(A, C, \re, \le)$
with $C \in \rema$, there is a bijective correspondence between
left $\cA$-linear maps $g: C \otimes C \to \cA$ and $k$-linear
maps $\hat{g}:C \otimes C \to A$ satisfying condition \eqref{cond}.
\end{lemma}
\begin{proof}
Note that, in view of the form of the left $A$-action in Corollary~\ref{cla}, the condition
\eqref{cond} is equivalent to
\begin{equation}\label{cond1}
\sum_\alpha a \sba \hat{g} (c \spa \otimes c') =
\hat{g}(c\sw 1 \ot c')\eps_C(c\sw 2 a).
\end{equation}
Given a $k$-linear map $\hat{g}$ satisfying condition (\ref{cond})
define $g: C \otimes C \to \cA$ as $g:= \lpr \circ (C \otimes
\hat{g}) \circ (\Delta_C \otimes C)$, so on elements $g(d \otimes
d')= \sum \sba d \sw 1 \otimes \hat{g}(d \sw 3 \otimes d') \sba
\eps_C ( d \sw 2 \spa)$. Using (\ref{WE2}),  (\ref{WE4}) and condition
(\ref{cond1}) we obtain, for all $d,d'\in C$,
$$
 \sum \sba
\hat{g}(d \sw 2 \otimes d') \sba
\eps_C  (d \sw 1 \spa)
= \sum \sba \eps_C  (d \sw 1 \spa) 1 \sba \hat{g} ( d \sw 2 \otimes d')
= \sum \sba 1 \sba \hat{g}(d \spa \otimes d')
= \hat{g} (d \otimes d'),
$$
hence
\begin{equation}\label{form}
g(d\ot d') = d\sw 1\ot \hat{g}(d\sw 2 \ot d').
\end{equation}
Next note that $\cA\coten{C}C\ot C$ consists of $k$-linear combinations of
$\suma c\sw 1\ot a\sba \ot c\sw 2\spa \ot d$, with $a\in A$ and $c,d\in C$, and compute
\begin{eqnarray*}
\suma g((c\sw 1\ot a\sba)\rhd ( c\sw 2\spa \ot d)) &=& \suma d(c\sw 1\ot \eps_C(c\sw 2
a\sba)\eps_C(c\sw 3\spa)d)\\
&=& g(c\sw 1\ot d)\eps_C(c\sw 2 a) = c\sw 1\ot\hat{g}(c\sw 2\ot d)\eps_C(c\sw 3 a)\\
&=& \suma c\sw 1\ot a\sba \hat{g}(c\sw 2\spa\ot d) = \suma (c\sw 1\ot a\sba )g( c\sw 2\spa \ot d),
\end{eqnarray*}
where the second equality follows by the fact that $C$ is a weak entwined module,
the third by \eqref{form}, the fourth by condition \eqref{cond1}. The final equality
is a consequence of \eqref{form} and the definition of product in $\cA$. This shows that
$g$ is a left $\cA$-module map.

For the converse, given a left $\cA$-linear map $g: C \otimes C \to
\cA$ define $\hat{g} : C \otimes C \to A$ to be $\hat{g}:= (\eps_C
\otimes A) \circ g$.  Observe that $\sum \sba d \sw 1 \otimes
a \sba \otimes d \sw 2 \spa \otimes d'$ lies in $\cA \cotc C
\otimes C$ for all $a\in A$ and $d \in C$. Apply
the map $\eps_C  \otimes A : \cA \to A$ to the
$\cA$-linearity condition of $g$
$$
\sum \sba (d \sw 1 \otimes a \sba) g (d \sw 2 \spa \otimes d') =
\sum \sba g ((d \sw 1 \otimes a \sba) \rhd (d \sw 2 \spa \otimes
d')), 
$$
and observe
that  $\eps_C  \otimes A $ is multiplicative with respect to the $C$-ring product in
$\cA,$ to conclude that $\hat{g}$
satisfies the required condition (\ref{cond}).

It remains to show that  the
given correspondence is one-to-one. Clearly, applying $\eps_C\ot A$ to
$g$ given in terms of $\hat{g}$ via equation \eqref{form},
one obtains back $\hat{g}$. On the other hand, since $g$ is left $C$-colinear,
$g = (C\ot \eps_C \ot A)\circ (C\ot g)\circ (\Delta_C \ot C)$, thus establishing the
converse correspondence.
\end{proof}

Using this lemma we are now able to prove Proposition~\ref{prop.co.ext}.

\begin{proof} (Proposition~\ref{prop.co.ext})  Suppose that there is a map
$\hat{g} :C\ot C\to A$ satisfying \eqref{cond} and \eqref{cotran}. By
Lemma~\ref{lemma2} there is
a corresponding left $\cA$-linear map $g: C\ot C\to\cA$, $c\ot c'\mapsto
c\sw 1\ot \hat{g}(c\sw 2\ot c')$. The condition \eqref{cotran} ensures that
$g$ is a retraction of $\beta$, hence
$ C \twoheadrightarrow B $ is a weak $A$-Galois coextension and
$C$ is  injective as a left $B$-comodule by Theorem~\ref{theorem1}.

Conversely, if $ C \twoheadrightarrow B $ is a weak $A$-Galois coextension
and
$C$ is  injective as a left $B$-comodule, then, by Theorem~\ref{theorem1}
there is a left $\cA$-module retraction $g$ of $\beta$. The map $\hat{g} =
(\eps_C\ot A)\circ g$ satisfies \eqref{cond} (by Lemma~\ref{lemma2}) and
\eqref{cotran} (since $g$ is a retraction of $\beta$).
\end{proof}

In the case where $\re$ is a bijective entwining structure
(non-weak!), $\psi_R$ is a bijective map (with the inverse
$\psi_L$), hence the condition \eqref{cond} means that $\hat{g}$
is left $A$-linear.
\begin{example}\label{homog}
Let $H$ be a Hopf algebra with bijective antipode $S$, and let $A$ be a right
$H$-coideal subalgebra of $H$, i.e.\ $A$ is a subalgebra of $H$ and
$\Delta_H(A) = A\ot H$. In this case $(A,H,\psi_R)$, with
$$
\psi_R : H\ot A\to A\ot H,\qquad h\ot a\mapsto a\sw 1\ot \ ha\sw 2,
$$
we have a bijective right entwining structure for which $H$ is an
entwined module. The inverse of $\psi_R$ is
$$
\psi_L: A\ot H\to H\ot A, \qquad a\ot h\mapsto hS^{-1}(a\sw 2)\ot \ a\sw 1,
$$
hence the induced left $A$-action on $H$ is $a\cdot h := h S^{-1}(a)$. Suppose that $A$
is a direct summand of $H$ as a left $A$-module (e.g.\ there is a strong connection in
$H$, cf.\ \cite[Theorem~2.5]{BrzHaj:Che}), and let $p: H\to A$ be a left $A$-linear retraction of $A\subseteq H$. Define
the map
$$
\hat{g}: H\ot H\to A, \qquad h\ot h'\mapsto p(S(h)h').
$$
Then the map $\hat{g}$ satisfies both \eqref{cond} and
\eqref{cotran}, hence $H \twoheadrightarrow B $ is an $A$-Galois
coextension and $H$ is  injective as a left $B$-comodule. In this
case $B=H/HA^+$, where $A^+ = A\cap \ker\eps_H$.
\end{example}

As a concrete illustration of Example~\ref{homog}, take $H = \mathcal{O}(SU_q(2))$,
the
algebra of (polynomial) functions on the quantum group $SU_q(2)$ \cite{Wor:twi}
  and $A=  \mathcal{O}(S^2_{q,s})$,  the
algebra of (polynomial) functions on the quantum two-sphere \cite{Pod:sph}.
$\mathcal{O}(SU_q(2))$ is known to be a coalgebra-Galois extension
of $\mathcal{O}(S^2_{q,s})$ with a strong connection
(explicitly constructed in \cite{BrzMaj:geo}). This implies that
 $\mathcal{O}(S^2_{q,s})$
is a direct summand
of $\mathcal{O}(SU_q(2))$ as a left $\mathcal{O}(S^2_{q,s})$-module.
The coinvariant coalgebra $B$ is spanned by countably many group-like
elements  (hence it can be identified with the
Hopf algebra $\mathcal{O}(S^1) = k[Z,Z^{-1}]$).
Consequently, $\mathcal{O}(SU_q(2))$ is an $\mathcal{O}(S^2_{q,s})$-Galois
coextension of $B$ and it is  injective as a $B$-comodule.

Proposition~\ref{prop.co.ext} can be used to characterise weak Galois
coextensions of self-injective algebras.
\begin{theorem}\label{KTS.thm}
Let $(A , C, \re, \le)$ be an invertible weak entwining structure
such that $C$ is a right entwined module,
and let $\cA$ be the $C$-ring corresponding to $(A, C, \re)$.
Suppose that the map $\beta: \cA \to C \otimes C$,
$c \otimes a \mapsto c \sw 1 \otimes
c \sw 2 a$ is injective.
 If $A$ is a left self-injective algebra, then
$ C \twoheadrightarrow B $ is a weak $A$-Galois coextension and
$C$ is  injective as a left $B$-comodule. Furthermore, if $A$ is a
separable algebra, then  $C$ is also $A$-equivariantly injective
as a left $B$-comodule (i.e., $C$ is an injective left $B$
comodule and the corresponding coaction has a retraction in $^B \M
_A$).
\end{theorem}

\begin{proof}
 Firstly  view $\cA$ as a left
$A$-module by
$$\xymatrix{A\ot \cA \ar[r]^{\le\ot A} & C\ot A\ot A \ar[r]^{C\ot \mu_A}&\cA.}
$$
This is easily seen to be well-defined, since $\cA = \im \lpr = \im p_L$.
Secondly, view $C \otimes C$ as a left $A$-module through the
composition
$$
A \otimes C \otimes C \xrightarrow{\le\ot C} \cA \otimes C
\xrightarrow{\overline{_{C \otimes C}\rho}} C \otimes C,
$$
i.e., use the left $A$-action in Corollary~\ref{cla}, $a\ot c\ot c'\mapsto ac\ot c'$.
Define the map
$$ r: \cA \to A, \qquad r(c \otimes a)=\eps_C (c)a,$$
and observe that, for all $b\in A$ and $c\ot a\in \cA$ (summation
suppressed for clarity),
$$
r(b  (c \otimes a)) = \sum_E \eps_C (c_E) b^E a = \sum_E  b \eps_C (c_E) 1^Ea
= b \eps_C (c)a.
$$
The second equality is by (\ref{le2}) and final equality since $c
\otimes a \in \cA$ implies that $\sum_E  c_E \otimes 1^E a=p_L(c
\otimes a)=c \otimes a$.  Hence $r \in \mathrm{Hom}_{A-}(\cA, A)$.
Next we prove that the map $\beta :\cA\to C\ot C$ is also $A$-linear. This is
done in a few steps. First, using \eqref{le3} note that, for all $a\in  A$ and $c\in C$,
\begin{equation}\label{eq*}
\sum_E {c\sw 1_E}\eps_C(a^Ec\sw 2) = ac.
\end{equation}
On the other hand, since $C\in \rem$ and $\psi_L\circ\psi_R =
\lpr$, we find that
\begin{eqnarray*}
\suma\eps_C(a\sba c\spa \sw 1)c\spa \sw 2 &=& \sum_{\alpha, E} \eps_C({c\spa \sw 1}_E {a\sba}^E)c\spa\sw 2 = \sum_{\alpha,\beta,E} \eps_C({c\sw 1\spb}_E{a_{\alpha\beta}}^E)c\sw 2\spa\\
&=& \sumab \eps_C(c\sw 1 a_{\alpha\beta})\eps_C(c\sw 2\spb)c\sw 3\spa = ca,
\end{eqnarray*}
where the second equality follows by the definition of the left
$A$-action in Corollary~\ref{cla}. We can combine this way of
expressing of right $A$-action on $C$ in terms of the left
$A$-action with the equality $\re\circ\le = \lpl$ and the fact
that $C\in\lem$, to find that, for all $c\in C$ and $a\in A$,
\begin{equation}\label{eq**}
\sum_E c_Ea^E = \sum_{\alpha, E} \eps_C({a^E}\sba{c_E}\spa\sw 1){c_E}\spa\sw 2
= \sum_E \eps_C(a^Ec\sw 2)\eps_C({c\sw 1}_E) = \eps_C(ac\sw 1)c\sw 2.
\end{equation}
Therefore, for all $a,b\in A$ and $c\in C$,
\begin{eqnarray*}
\beta(bp_L(c\ot a)) &=& \sum_E {c_E}\sw 1\ot {c_E}\sw 2b^E a = \sum_{E,F}
{c\sw 1}_E\ot {c\sw 2}_F b^{EF}a\\
&=& \sum_E {c\sw 1}_E\eps_C(b^Ec\sw 2)\ot c\sw 3a = bc\sw 1\ot c\sw 2a =
b\beta(p_L(c\ot a)),
\end{eqnarray*}
where the second equality is by \eqref{le3}, the third by
\eqref{eq**} and the fourth by \eqref{eq*}. This proves that
$\beta$ is a left $A$-linear map, and thus, in view of the
self-injectivity of $A$, we are led to
an exact sequence
$$ \mathrm{Hom}_{A-}(C \otimes C, A) \xrightarrow{\beta^*} \mathrm{Hom}_{A-}(\cA,
A)\rightarrow 0$$ and so there exists $\hat{g} \in
\mathrm{Hom}_{A-}(C \otimes C, A)$ s.t. $\beta ^* \circ \hat{g} =
\hat{g} \circ \beta = r$.  By construction, $\hat{g}$ satisfies
condition \eqref{cotran} and it is left
$A$-linear, hence   \eqref{cond} holds. By Proposition~~\ref{prop.co.ext},
$ C \twoheadrightarrow B $ is a weak $A$-Galois coextension and
$C$ is  injective as a left $B$-comodule.

Now suppose furthermore that $A$ is a separable algebra and let $e
= e_1 \otimes e_2 \in A \otimes A$ denote the separability element
(summation suppressed). To show that $C$ is $A$-equivariantly
injective as a left $B$-module we need to show that there exists a
retraction of the left $B$-coaction, given in Corollary
~\ref{proof1}, in $ ^B \M _A$.  The injectivity of $C$ as a left
$B$-module guarantees that there is a left $B$-colinear map $\hat
{\lambda} : B \otimes C \to C$ such that $\hat{\lambda} \circ
{^C\! \varrho} = C$. From this we can construct
$$
 \lambda : B \otimes C \to C, \qquad \lambda = {\roC} \circ
(\hat{\lambda} \otimes A) \circ (B \otimes {\roC} \otimes A) \circ
(B \otimes C \otimes e).
$$
Now observe that ${ \roC}:C \otimes A \to C$ is a left
$B$-comodule map because $^C\! \varrho : C \to B \otimes C$, given
in Corollary~\ref{proof1}, is right $\cA$-linear and the
correspondence given in the third part of
Theorem~\ref{thm.weak.C-ring} allows the right $A$-action on $C$
to be viewed as some right $\cA$-action.
  With this in mind it is clear that $\sigma$ is left
$B$-colinear, since it is a composition of $B$-colinear maps.
That it is a right $A$-linear map follows by the fact that
$ea=ae$, for all $a\in A$. It only remains to show that this map
is indeed a retraction for the left $B$-coaction. Just compute
\begin{eqnarray*}
\lambda \circ {^C \rho (c)} &=& \hat{\lambda}(c \zz {-1} \otimes c
\zz 0  e_1)  e_2\\
&=& \hat{\lambda} ( (c  e_1) \zz {-1} \otimes ( c  e_1)
\zz 0)  e_2\\
&=& c  e_1 e_2 = c,
\end{eqnarray*}
where the second equality follows from the left $B$-colinearity of
the right $A$-action, the third because $\hat{\lambda}$ was chosen
to be a splitting of the coaction and the final equality from the
properties of the separability element.
\end{proof}

Theorem~\ref{KTS.thm} is  a dual version of \cite[Theorems 5.1, 6.1]{BrzTur:str},
thus a dualisation of each in the long chain of generalisations of
the Kreimer-Takeuchi theorem \cite[Theorem~1.7]{KreTak:hop}
for Hopf-Galois extensions. In particular, in its self-injective part, the
non-weak case corresponds to \cite[Theorem~3.1]{SchSch:gen}, the proof of which
lends the idea for the proof of Theorem~\ref{KTS.thm}. Since any
quasi-Frobenius algebra is self-injective, Theorem~\ref{KTS.thm} implies also
a dual version of
\cite[Theorem~3.1]{BeaDas:Gal}. In particular, this is applicable to
extensions of finite dimensional weak Hopf algebras. Any such weak Hopf algebra
$H$ has a bijective antipode by \cite[Theorem~2.10]{Boh:wea} thus the
weak entwining structure $(H,C,\re)$ corresponding to a right $H$-module
coalgebra $C$ is invertible. Furthermore, a finite dimensional weak
Hopf algebra is quasi-Frobenius by \cite[Theorem~3.11]{Boh:wea}. Hence
Theorem~\ref{KTS.thm} implies that the injectivity of the canonical map
$\beta$ is sufficient for a coextension $C$ of a finite dimensional weak Hopf algebra
$H$ to be a weak Hopf-Galois coextension.

\end{document}